\documentclass{article}

\usepackage{todonotes}

\usepackage[utf8]{inputenc}
\usepackage{natbib}
\input{rose.sty}

\usepackage{algorithm}
\usepackage{algpseudocode}
\usepackage{xcolor}
\usepackage{enumitem}
\usepackage{todonotes}

\newtheorem*{corollary*}{Corollary}
\newtheorem*{lemma*}{Lemma}

\newcommand{\hamilton}{\mathcal{H}}

\newcommand{\hamiltonGCW}{\mathcal{H}^{\GCW}_{\beta}}
\newcommand{\polyGCW}{h_{\beta}^\GCW}
\newcommand{\rateGCW}{L^\GCW_\beta}
\newcommand{\ratemaximaGCW}{M^\GCW_\beta}

\newcommand{\phasemuGCW}[1]{\mathcal{M}_{\beta,#1}^{\GCW}}

\title{Approximate FKG inequalities for phase-bound spin systems,
with applications to central limit theorems for exponential random graphs} 
\author{
    Satyaki Mukherjee \\
    National University of Singapore \\
    \texttt{satyaki.mukhyo@gmail.com}
    \and
    Vilas Winstein \\
    University of California, Berkeley \\
    \texttt{vilas@berkeley.edu}
}

\begin{document}

\maketitle

\begin{abstract}
The Fortuin--Kasteleyn--Ginibre (FKG) inequality is an invaluable tool in monotone spin systems satisfying the FKG lattice condition,
which provides positive correlations for all coordinate-wise increasing functions of spins.
This inequality has numerous applications and plays an integral role in the proof of various central limit theorems (CLTs),
including recent work on ferromagnetic exponential random graph models (ERGMs) wherein a Hamiltonian tilt promotes
the presence of small subgraphs like triangles \cite{ganguly2024sub,fang2025normal,winstein2025quantitative}.
However, the FKG lattice condition fails to hold when confining a spin system to a particular
phase in the low-temperature regime of parameters.
Thus it is not a priori clear if each phase internally has positive correlations for increasing functions,
or if the positive correlations in the overall model (which is a mixture of phases) arise primarily from the global choice of phase.

In this article, we show that the individual phases in ERGMs do indeed satisfy an approximate form of the FKG inequality internally.
We use this to finish the proof of various CLTs within each individual phase in the phase-coexistence regime,
answering a question posed by Bianchi, Collet, and Magnanini \cite{bianchi2024limit}.
We present the FKG inequality for ERGMs as a consequence of a more general result which holds under certain inputs related to
\emph{metastable mixing}; we expect this general result to be widely applicable, and we devote a section to spelling out the details
of its application to a class of generalized higher-order ferromagnetic Curie--Weiss models where the necessary inputs are relatively transparent.
\end{abstract}

\setcounter{tocdepth}{2}
\tableofcontents

\section{Introduction}
\label{sec:intro}

The FKG inequality is an essential tool in monotone spin systems which allows one to deduce nonnegative correlations for all
increasing functions from a simpler positive-correlation condition.
Abstractly, if $\Omega = \calA^N$ for some ordered set $\calA$ and integer $N$, and if $\mu$ is a probability measure on $\Omega$ satisfying
\begin{equation}
\label{eq:lattice}
    \mu(x \wedge y) \mu(x \vee y) \geq \mu(x) \mu(y),
\end{equation}
then for any $f, g : \Omega \to \R$ which are increasing coordinate-wise, we have $\Cov_\mu[f,g] \geq 0$.
Here $x \wedge y$ and $x \vee y$ denote the coordinate-wise minimum and maximum respectively, for $x, y \in \Omega$.
Initially studied by Harris \cite{harris1960lower} as well as Fortuin, Kasteleyn, and Ginibre \cite{fortuin1971correlation},
the FKG inequality allows for the construction of infinite-volume Gibbs measures for the lattice Ising model \cite{dobrushin1968gibbsian}, 
is a key input for the Russo--Seymour--Welsh theory which gives bounds on critical probabilities in percolation \cite{russo1978note,seymour1978percolation},
can be used to derive central limit theorems in various contexts \cite{newman1980normal,bulinski1998asymptotical}, and has many other applications.

Many examples of spin systems on $\Omega$ satisfy the FKG condition \eqref{eq:lattice} in their so-called \emph{ferromagnetic} phases,
where local interactions promote alignment of spins.
Such models are relatively simple at high temperatures (weak interaction) and are governed primarily by entropy,
leading to approximate independence of spins.
However, at low temperatures they exhibit more complex behavior which is a subject of modern interest in statistical mechanics,
and a common theme has emerged whereby low-temperature spin systems may often be separated as mixtures of
\emph{phases}, each of which has different properties, which behave roughly like high-temperature systems 
when considered separately from each other.
In addition, since the phases are often \emph{metastable} under the dynamics, meaning phase transitions take a very long time,
it makes physical sense to consider each phase as a model in its own right.
These ideas have been made precise in a variety of situations including mean-field and lattice Ising models
\cite{samanta2024mixing,gheissari2022low}, quantum error-correcting codes \cite{bergamaschi2025rapid,bergamaschi2025structural},
some disordered models such as spin glasses \cite{gamarnik2021overlap,sellke2024threshold,huang2025weak}, and ferromagnetic exponential random graph models (ERGMs)
\cite{bresler2024metastable,winstein2025concentration,winstein2025quantitative} which are the focus of this article and will be
discussed in Section \ref{sec:intro_setup_ergm} below.
These works use the study of separate phases to give improved sampling algorithms via sampling in stages, explain why other sampling
algorithms fail, or give other deep insights into the structure of these models.

Considering the phases of low-temperature spin systems separately poses an interesting mathematical challenge, especially
because many of the tools used at high temperatures break down.
As a first example to demonstrate the problem,
let us consider the Curie--Weiss model, where $\Omega = \{ -1, +1 \}^N$ and spins interact via a quadratic potential, i.e.\
\begin{equation}
    \mu(x) \propto \Exp{- \frac{1}{2T} \sum_{i,j=1}^N |x(i) - x(j)|^2},
\end{equation}
which encourages alignment of spins with an interaction strength which is inversely proportional to temperature $T$.
At low temperatures, the mass of the measure is split between two phases with magnetization (average spin)
close to either $+ m_*$ or $- m_*$, for some $m_* \in (0,1)$.
The measure $\mu$ satisfies the FKG condition \eqref{eq:lattice} globally, but if we consider a single phase measure, i.e.\ 
condition $\mu$ to have magnetization close to $+ m_*$, say, then \eqref{eq:lattice} will no longer hold as taking the coordinate-wise
maximum or minimum may drastically change the magnetization, so the left-hand side of \eqref{eq:lattice}
may be zero even if the right-hand side is positive.
As such, it has remained an interesting question whether the positive correlations in the low-temperature Curie--Weiss model
arise from local interactions within the phases, or from the overall choice of phase, i.e.\ positive or negative magnetization.
In the present article, we show that the first statement holds by exhibiting (almost) positive correlations of increasing functions
for the phase measures of a class of generalized Curie--Weiss models, including the standard Curie--Weiss model as a special case.

We also apply our techniques to the class of ferromagnetic ERGMs mentioned above, which are our main focus in the present article.
In these models, $\Omega = \{0,1\}^{\binom{n}{2}}$ is thought of as the set of labeled graphs on $n$ vertices,
the spins are the edges of the graph, and they interact through the presence of small subgraphs such as triangles.
As mentioned, a full definition will be provided in Section \ref{sec:intro_setup_ergm} below, but for now note that
much like the Curie--Weiss model there is a low-temperature parameter regime where the mass of the measure may be split between
two (or more) phases.
As with the generalized Curie--Weiss models,
we show that within each phase an ERGM still exhibits (almost) positive correlations of increasing functions.

For ERGMs, the FKG inequality has been an essential tool in proving central limit theorems for the edge count and other observables
\cite{ganguly2024sub,fang2025normal,winstein2025quantitative}, and as an application of our methods, in Theorem~\ref{thm:main_clt} we are able to extend such results to
the phase measures in cases which have eluded analysis thus far, answering a question of \cite{bianchi2024limit}
for the edge--triangle model in the phase-coexistence regime.
We mention briefly that such central limit theorems have a wide array of applications which include parameter estimation for ERGMs \cite{mukherjee2013statistics},
a result which should prove useful for practitioners working with these models in the phase-coexistence regime.
However, to keep the scope of this paper focused, we do not pursue this presently.

Finally, we mention that our techniques are quite flexible, and we have consolidated them into a general result
(Theorem \ref{thm:general} below), which we apply to both models mentioned above.
We expect this general result to be applicable in various other contexts of interest beyond the two covered in the present article;
in fact, one may view our application to the generalized Curie--Weiss models as a template,
laying out in detail which inputs are necessary when applying our general result.

\subsection{Problem setup}
\label{sec:intro_setup}

Let us now properly define the two models we will consider in this article and introduce the
decomposition of these models into their phase measures as discussed above.
For the generalized Curie--Weiss models (GCWMs), we give a complete description in the next section,
but for the exponential random graph models (ERGMs) a fully rigorous understanding will have to wait
until Section \ref{sec:ergm_background_ldp}, although we give the basic idea in
Section \ref{sec:intro_setup_ergm} below.

\subsubsection{Generalized Curie--Weiss models}
\label{sec:intro_setup_gcwm}

A generalized Curie--Weiss model is parameterized by a polynomial
\begin{equation} \label{eq:polyGCW}
    \polyGCW(m) = \sum_{j=1}^K \beta_j m^j,
\end{equation}
or in other words by a fixed $K$ and a parameter vector $\beta = (\beta_1, \dotsc, \beta_K) \in \bbR^K$.
Given this polynomial, we define the \emph{Hamiltonian function} on $x \in \{0,1\}^N$ by setting
\begin{equation}
\label{eq:hamiltonGCW}
    \hamiltonGCW(x) \coloneqq \polyGCW(m(x)),
\end{equation}
where $m(x)$ is the \emph{mean magnetization}, given by
\begin{equation}
\label{eq:magnetization}
    m(x) = \frac{1}{N} \sum_{i=1}^N x(i).
\end{equation}
With the Hamiltonian \eqref{eq:hamiltonGCW}, we may define
a probability distribution on $x \in \{ 0, 1 \}^N$ by 
\begin{equation}
\label{eq:cw_full}
    \nu^\GCW_\beta(x) \propto \Exp{N \cdot \hamiltonGCW(x)}.
\end{equation}
The scaling by $N$ ensures that the measure $\nu^\GCW_\beta$ is nontrivial in the sense that
it is neither approximately uniform nor concentrated on the maximizers of $\hamiltonGCW$.
In the present article we make the additional assumption that $\beta_j \geq 0$ for each $j \geq 2$.
This is the \emph{ferromagnetic} assumption, which ensures positive correlations and an FKG
inequality for the full measure $\nu^\GCW_\beta$.
Note that we still allow $\beta_1$ to be negative.

Note also that the standard Curie--Weiss model is obtained, after making the linear change of variables
to $y = 2x-1 \in \{-1,+1\}^N$ by taking $K=2$ and setting
the external magnetic field to $\frac{\beta_1 + \beta_2}{2}$ and the inverse temperature
to $\frac{\beta_2}{4}$.
Note that we prefer to work with $\{0,1\}^N$ instead of $\{-1,+1\}^N$ in the present article simply for
consistency with the ERGMs which will be defined in the next section.

Just as the standard Curie--Weiss model satisfies a large deviations principle which shows that the magnetization
concentrates around some finite set of values, the same behavior can be observed in generalized Curie--Weiss models,
as will be shown in Proposition \ref{prop:gcw_ldp} below.
The possible magnetizations in the large deviations principle are the maximizers $m_*$ of the following rate function:
\begin{equation}\label{eq:rateGCW}
    \rateGCW(m) = \polyGCW(m) - m \log m - (1-m) \log(1-m).
\end{equation}
Letting $M_\beta^\GCW$ denote the set of these maximizers, Proposition \ref{prop:gcw_ldp} states that if
$\tilde{X} \sim \nu_\beta^\GCW$ then
for any $\eta > 0$ there are constants $C(\eta),c(\eta) > 0$ such that
\begin{equation}
\label{eq:cw_ldp_initial}
    \P\left[\inf_{m_* \in M_\beta^\GCW} |m(\tilde{X}) - m_*| > \eta \right] \leq C(\eta) e^{-c(\eta) N}.
\end{equation}
We thus find that the measure $\nu^\GCW_\beta$ decomposes
as a mixture of \emph{phase measures} with magnetizations near $m_* \in M_\beta^\GCW$, up to total-variation error $e^{-\Omega(N)}$.
To be precise, for any $m_* \in M_\beta^\GCW$, we define the phase measure $\calM_{\beta,m_*}^\GCW$ to be the measure
$\nu^\GCW_\beta$ conditioned on the set
\begin{equation}
\label{eq:gcw_band}
    \{ x : |m(x) - m_*| \leq \eta \}.
\end{equation}
Thus these measures depend implicitly on a choice of $\eta>0$, but we will suppress this from our notation as our results
do not depend on $\eta$ as long as it is small enough that the sets \eqref{eq:gcw_band} are disjoint for $m_* \in M_\beta^\GCW$.

These phase measures, considered separately from one another, are our main object of study, and Theorem \ref{thm:main_gcwm}
below states that each one satisfies an approximate version of the FKG inequality.
For technical reasons, we must restrict our attention to the so-called \emph{non-critical} case,
where $L^\GCW_\beta$ is \emph{strictly concave} at the global maximizer $m_*$; in other words, we define
\begin{equation}\label{eq:ratecritGCW}
    U^\GCW_\beta = \left\{ m_* \in \ratemaximaGCW : \left. \frac{d^2}{dm^2} L^\GCW_\beta(m) \right|_{m=m_*} < 0 \right\},
\end{equation}
and in the sequel we will exclusively consider $\phasemuGCW{m_*}$ for $m_* \in U^\GCW_\beta$.

\subsubsection{Exponential random graph models}
\label{sec:intro_setup_ergm}

An exponential random graph model is a probability distribution on $n$-vertex labeled graphs $x$,
which we identify with elements of $\{0,1\}^N$ for $N = \binom{n}{2}$.
We index this set not by integers but rather by two-element subsets of $[n] = \{1,\dotsc,n\}$,
which we identify with the edges of the complete graph $K_n$.
We denote this set of edges by $\edgeset$, and for any $e \in \edgeset$ we think of $x(e)$ as the indicator that the edge $e$
is in the graph $x$.
For any fixed finite graph $G$, let us denote by $N_G(x)$ the number of
\emph{homomorphisms} of $G$ in $x$, i.e.\ the number of maps $\calV(G) \to [n]$ which map
edges of $G$ to edges of $x$.
Here we use $\calV(G)$ and $\calE(G)$ to denote the vertex and edge set of a graph respectively.
Note that the maps counted by $N_G(x)$ \emph{do not} need to map non-edges to non-edges, i.e.\ they are not
necessarily \emph{induced} homomorphisms.
Additionally, we do not require them to be injective, but this is a relatively mild point.
Let us further define
\begin{equation}
    t(G,x) = \frac{N_G(x)}{n^{|\calV(G)|}},
\end{equation}
which is the \emph{homomorphism density} of $G$ in $x$, or in other words the probability that a randomly
chosen map $\calV(G) \to [n]$ is a homomorphism.

For any sequence of fixed finite graphs $G_0, \dotsc, G_K$, where $G_0$ is always assumed to
be a single edge, and any parameter vector
$\beta = (\beta_0, \dotsc, \beta_K) \in \R^{1+K}$, we define the ERGM Hamiltonian as follows:
\begin{equation}
\label{eq:ERGM_hamiltonian}
    \hamiltonERG(x) = \sum_{j=0}^K \beta_j \cdot t(G_j,x),
\end{equation}
and define the full ERGM distribution as
\begin{equation}
\label{eq:full_ergm_measure}
    \nu^\ERG_\beta(x) \propto \Exp{n^2 \cdot \hamilton^\ERG_\beta(x)}.
\end{equation}
As in the case of GCWMs, the scaling by $n^2$ ensures nontriviality.
We again make a \emph{ferromagnetic} assumption, that $\beta_j \geq 0$ for all $j \geq 1$,
but we do still allow the coefficient $\beta_0$ of the single edge homomorphism density $t(G_0,x)$ to be negative.

Exponential random graph models also satisfy a large deviations principle resulting in a
phase decomposition, although it is somewhat more difficult to state precisely.
For now, let us just note that under the ferromagnetic assumption,
the phases or metastable wells correspond to maximizers $p_*$ of
the function
\begin{equation}
\label{eq:Lerg_def}
    L^\ERG_\beta(p) = h_\beta^\ERG(p) - \frac{1}{2} \left(
    p \log p + (1-p) \log(1-p) \right).
\end{equation}
Here we define the polynomial $h_\beta^\ERG$ by
\begin{equation}
    h_\beta^\ERG(p) = \sum_{j=0}^K \beta_j p^{|\calE_j|},
\end{equation}
where $\calE_j$ denotes the edge set of $G_j$.
This arises from the natural interpretation of the subgraph density 
for the \emph{density-$p$ constant graphon} $W_p$; the notion
of graphons will be discussed in Section \ref{sec:ergm_background_ldp}, but for now this constant
graphon should be thought of as a limit of Erd\H{o}s--R\'enyi graphs under the \emph{cut distance},
$\db$ which will also be discussed in Section \ref{sec:ergm_background_ldp}.
Letting $M^\ERG_\beta \subseteq [0,1]$ denote the set of maximizers of $L^\ERG_\beta$, the large
deviations principle for exponential random graph models (due to \cite{chatterjee2013estimating})
states that
if $\tilde{X} \sim \nu^\ERG_\beta$, then
for any $\eta > 0$ there are constants $C(\eta), c(\eta) > 0$ for which
\begin{equation}
    \P \left[ \inf_{p_* \in M^\ERG_\beta} \db(\tilde{X},W_{p_*}) > \eta \right]
    \leq C(\eta) e^{-c(\eta) n^2}.
\end{equation}
Thus we may define the phase measure or metastable well $\calM^\ERG_{\beta,p_*  }$, for each
$p_* \in M^\ERG_\beta$, as the measure $\nu^\ERG_\beta$ conditioned on the set
\begin{equation}
    \left\{ x : \db(x,W_{p_*}) \leq \eta \right\},
\end{equation}
for some $\eta > 0$ chosen small enough so that these sets are disjoint as $p_*$ varies in
$M^\ERG_\beta$.
Again, the details of this construction will be made precise in Section \ref{sec:ergm_background_ldp} below.

As with the GCWMs, we consider these phase measures separately, and Theorem \ref{thm:main_ergm} below states that
each one satisfies an approximate version of the FKG inequality.
Note that as in the case of GCWMs, here we also only consider \emph{non-critical}
phases, where $L^\ERG_\beta$ is strictly concave, and we define $U^\ERG_\beta \subseteq M^\ERG_\beta$
as the set of maximizers $p_*$ where the second derivative of $L^\ERG_\beta$ is strictly negative.
\subsection{Results}
\label{sec:intro_results}

We now state our main results, which give almost nonnegative correlations for coordinate-wise increasing functions
of states under the phase measures of generalized Curie--Weiss models and exponential random graph models.
Both Theorems \ref{thm:main_gcwm} and \ref{thm:main_ergm} below follow from the more general
Theorem \ref{thm:general}, but the statement of that result is rather long and relegated
to the next section instead of being presented here.
In both results, $\| f \|_\infty$ denotes the maximum possible absolute value that the function $f$ can take.

First, Theorem \ref{thm:main_gcwm} covers the case of GCWMs.

\begin{theorem}
\label{thm:main_gcwm}
Let $\phasemuGCW{m_*}$ denote a phase measure in an $N$-spin generalized Curie--Weiss model,
for some $m_* \in U^\GCW_\beta$, recalling these definitions from Section \ref{sec:intro_setup_gcwm},
and consider a random spin configuration $X \sim \phasemuGCW{m_*}$.
There is some constant $c > 0$ such that for any coordinate-wise increasing $f$ and $g$,
\begin{equation}
    \Cov[f(X),g(X)] \geq - e^{-c N} \| f \|_\infty \| g \|_\infty.
\end{equation}
\end{theorem}

Next, Theorem \ref{thm:main_ergm} covers the case of ERGMs; we remark again that the proper rigorous
definition of the phase measures in this case will be given in Section \ref{sec:ergm_background_ldp}.

\begin{theorem}
\label{thm:main_ergm}
Let $\calM^\ERG_{\beta,p_*}$ denote a phase measure in an $n$-vertex exponential random graph model,
for some $p_* \in U^\ERG_\beta$, recalling these definitions from Section \ref{sec:intro_setup_ergm},
and consider a random graph $X \sim \calM^\ERG_{\beta,p_*}$.
There is some constant $c > 0$ such that for any coordinate-wise increasing $f$ and $g$,
\begin{equation}
    \Cov[f(X),g(X)] \geq - e^{-c n} \| f \|_\infty \| g \|_\infty.
\end{equation}
\end{theorem}

Note that the exponentially small error factor in the two theorems is actually different; in the
case of GCWMs we get an error which is exponentially small in the number of spins $N$, but for ERGMs
the error is only exponentially small in the number of vertices $n$, whereas the number of ``spins'' 
or potential edges is $N = \binom{n}{2}$.
This is ultimately due to the existence of small metastable regions of size $e^{-\Theta(n)}$
within the phase regions of ERGMs; see \cite[Theorem 3.3]{bresler2024metastable} for an explanation
of this phenomenon.

In the case of ERGMs, the positive correlation of increasing functions has been a useful tool for
deriving quantitative central limit theorems for various observables such as the total edge count, degree of a
particular vertex, or global/local subgraph counts \cite{fang2025normal,winstein2025quantitative}.
As such, these theorems have so far been restricted to situations where the phase measures are
close to the full measure, which happens only when there is exactly one relevant phase measure, i.e.\
when $|U_\beta^\ERG| = |M_\beta^\ERG| = 1$, which we term \emph{phase uniqueness}.
In \cite{winstein2025concentration}, a method was found to get around this restriction in certain
cases, but an important case which was left out was the edge--triangle or Strauss model, where,
in the notation of Section \ref{sec:intro_setup_ergm}, $K = 1$ and $G_1$ is a triangle.

The edge--triangle model is perhaps the most well-studied ERGM aside from the two-star ERGM (where $K=1$ and $G_1$ is a two-star or wedge graph)
which behaves similarly to an Ising model \cite{mukherjee2013statistics}.
Indeed a central limit theorem for the edge count in the edge-triangle model was
already shown in \cite{bianchi2024limit} in the phase uniqueness case, albeit with different methods
to \cite{fang2025normal,winstein2025quantitative} which do not use the FKG inequality and also do not
give quantitative bounds on the rate of convergence. 
In \cite{bianchi2024limit} it was asked whether the phase coexistence case could be handled, and our Theorem \ref{thm:main_ergm}
provides the key input needed to adapt the methods of \cite{fang2025normal,winstein2025quantitative} to this setting.
Thus we may state an unconditional central limit theorem for ERGM phase measures.

\begin{theorem}
\label{thm:main_clt}
Let $\calM^\ERG_{\beta,p_*}$ denote a phase measure in an $n$-vertex exponential random graph model,
for some $p_* \in U^\ERG_\beta$, and consider a random graph $X \sim \calM^\ERG_{\beta,p_*}$.
There is an explicit deterministic quantity $\sigma_n$ depending only on the ERGM specification
and the choice of $p_*$ for which
\begin{equation}
    \distance \left( \frac{|\calE(X)| - \E[|\calE(X)|]}{\sigma_n}, Z \right) \lesssim n^{-\frac{1}{2} + \eps}.
\end{equation}
Here $Z$ is a standard normal random variable, $|\calE(X)|$ denotes the number of edges in $X$,
and $\distance$ may be either the 1-Wasserstein distance or the Kolmogorov distance between $\bbR$-valued
random variables.
\end{theorem}

The deterministic quantity $\sigma_n$ will be explicitly given in Section \ref{sec:covariance_clt}.
In that section, we also state (or at least mention) a variety of other central limit theorems for quantities such as the degree of
a particular vertex, the overall count of a particular subgraph, and the count of subgraphs at a particular vertex.
All of these results follow using the same machinery developed in \cite{fang2025normal,winstein2025quantitative},
now using Theorem \ref{thm:main_ergm} as a crucial input to handle the phase-coexistence case.

\begin{remark}
As an immediate statistical application of Theorem \ref{thm:main_clt}, one may derive
a corresponding central limit theorem for the estimator $\hat{\theta}$ of the parameter
in any one-parameter family of parameter vectors $\beta(\theta)$,
as long as the curve $\beta : \R \to \R^{1+K}$ \emph{does not}
cross the phase transition surface for the ERGM in question,
and the derivative $\frac{d}{d\theta} p_*(\beta(\theta))$ is nonzero.
The novelty here is that the CLT in the phase-coexistence regime allows one to
consider curves $\beta$ which \emph{do enter} the phase transition surface, as long
as they do not cross from one side to the other.
\end{remark}

\subsection{Proof ideas}
\label{sec:intro_iop}

\begin{figure}
    \centering
    \begin{tikzpicture}[scale=1, baseline={(0,0)}]
        \fill[gray!20] (-3,-1) rectangle (3,1);
        
        \draw[gray] (-3,1) -- (3,1);
        \draw[gray] (-3,-1) -- (3,-1);
        
        \coordinate (top) at (0.25,2.25);
        \coordinate (right) at (2,0.5);
        \coordinate (bottom) at (-0.25,-1.75);
        \coordinate (left) at (-2,0);
        
        \draw[thick] (top) -- (right) -- (bottom) -- (left) -- cycle;
        
        \fill (top) circle (2pt);
        \fill (bottom) circle (2pt);
        \fill (right) circle (2pt);
        \fill (left) circle (2pt);
        
        \node[above] at (top) {$x \vee y$};
        \node[below] at (bottom) {$x \wedge y$};
        \node[left] at (left) {$x$};
        \node[right] at (right) {$y$};
    \end{tikzpicture}
    \hspace{2cm}
    \begin{tikzpicture}[scale=1, baseline={(0,0)}]
        \fill[gray!20] (-3,-1) rectangle (3,1);
        
        \draw[gray] (-3,1) -- (3,1);
        \draw[gray] (-3,-1) -- (3,-1);
        \draw[dashed,gray] (-3,0.75) -- (3,0.75);
        \draw[dashed,gray] (-3,-0.75) -- (3,-0.75);
        
        \coordinate (top) at (0.5,0.9);
        \coordinate (right) at (0.9,0.5);
        \coordinate (bottom) at (-0.5,-0.9);
        \coordinate (left) at (-0.9,-0.5);
        
        \draw[thick] (top) -- (right) -- (bottom) -- (left) -- cycle;
        
        \fill (top) circle (2pt);
        \fill (bottom) circle (2pt);
        \fill (right) circle (2pt);
        \fill (left) circle (2pt);
        
        \node[above] at (top) {$x \vee y$};
        \node[below, yshift=-0.5] at (bottom) {$x \wedge y$};
        \node[left] at (left) {$x$};
        \node[right] at (right) {$y$};
    \end{tikzpicture}
    \caption{
    Left: even if $x$ and $y$ are in the metastable well, $x \wedge y$ or $x \vee y$
    may not be.
    Right: to remedy this, we assume that $x$ and $y$ are \emph{almost} ordered, which
    means that $x \wedge y$ is close to either $x$ or $y$, and similarly for $x \vee y$.
    To make this argument work, we also need to restrict $x$ and $y$ to be a bit away from
    the boundary of the metastable well, so that nearby points remain in the well.
    Note that in these pictures, we visualize the relation $x \leq y$ as having $x$ lie
    within the cone under $y$ with sides at 45-degree angles from the horizontal.
    }
    \label{fig:diamonds}
\end{figure}
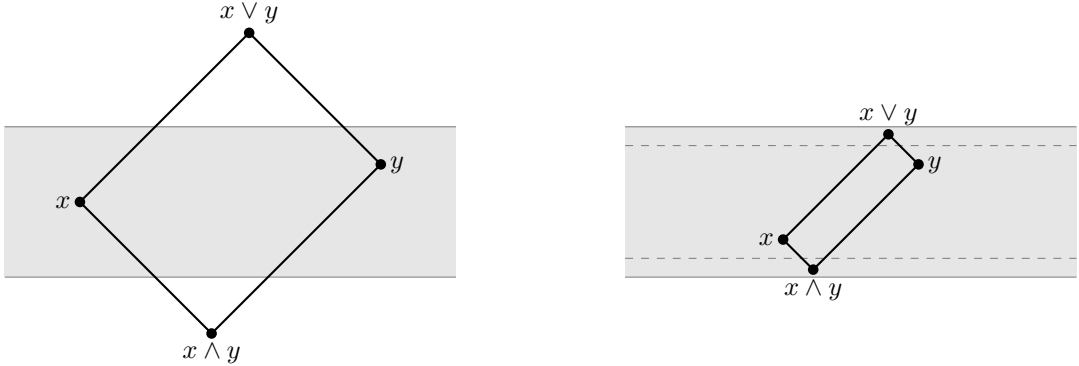

We now outline our arguments leading to Theorems \ref{thm:main_gcwm} and \ref{thm:main_ergm}.
Let us consider a measure $\mu$ on the Boolean hypercube $\{0,1\}^N$,
and recall that the FKG inequality implies positive correlations of increasing functions if
\begin{equation}
    \mu(x \wedge y) \mu(x \vee y) \geq \mu(x) \mu(y),
\end{equation}
which is the so-called FKG condition previously mentioned in \eqref{eq:lattice} above.
Now, even if $\mu$ satisfies this condition, if we condition the measure on a subset
of $\{ 0, 1 \}^N$, the condition may immediately be destroyed.
Indeed, let us consider the case of the generalized Curie--Weiss models introduced in
Section \ref{sec:intro_setup} above.
A metastable well in such a model is simply a ``band'' in the hypercube of the form
\eqref{eq:gcw_band}.
But if $x$ and $y$ lie in this band, then it is quite possible that $x \wedge y$ or $x \vee y$
do not, meaning that the left-hand side of the FKG condition is $0$ under the conditioned
measure while the right-hand side is not (see Figure \ref{fig:diamonds}, left).

To remedy this, we would like to consider only $x$ and $y$ such that $x \wedge y$ and $x \vee y$
remain in the band, and one simple way to achieve this is to work under the assumption
that $x \wedge y$ and $x \vee y$ are both close to either $x$ or $y$ themselves.
One way to achieve this is to assume that $x$ and $y$ are almost \emph{ordered}, i.e.\
that we almost have $x \leq y$ or $x \geq y$ in the partial order of the hypercube.
If we also restrict $x$ and $y$ to lie slightly away from the boundary of the metastable well,
then $x \wedge y$ and $x \vee y$ cannot escape it (see Figure \ref{fig:diamonds}, right).

To get away with this weaker and more ``local'' form of the FKG condition, we will also require
that the Glauber dynamics within the metastable well (to be described shortly) exhibits some form of \emph{contraction}, which leads
to rapid mixing.
To see why this is relevant, note that one proof of the standard FKG inequality due to Holley \cite{holley1974remarks}
goes via the analysis of Glauber dynamics to show that a tilted version of $\mu$ stochastically
dominates the original measure.
Let us next give a very brief outline of this argument; we refer the reader to
\cite[Section 3.10.3]{friedli2017statistical} for a more complete exposition as well.

Recall that in the standard FKG inequality we would like to show that $\Cov_\mu[f,g] \geq 0$
for increasing functions $f$ and $g$.
Let us consider the tilted measure given by
\begin{equation}
    \tilde{\mu}(x) \coloneqq \frac{g(x)}{\E_\mu[g]} \mu(x)
\end{equation}
(we may assume that $g$ is positive by shifting it by a constant, which does not change the covariance).
Then we have
\begin{equation}
    \Cov_\mu[f,g] = \E_\mu[g] \left( \E_{\tilde{\mu}} [f] - \E_\mu[f] \right),
\end{equation}
and so, since $f$ is increasing, it suffices to show that $\tilde{\mu}$ stochastically dominates
$\mu$, in the sense that $\tilde{X} \sim \tilde{\mu}$ and $X \sim \mu$ may be coupled so that
$\tilde{X} \geq X$ almost surely.

To construct such a coupling, we run Glauber dynamics, where a uniformly random coordinate of
the state is resampled (with the correct conditional distribution) at each time step.
Let us denote the stationary Glauber dynamics with respect to $\mu$ by $(X_t)$, and with respect
to $\tilde{\mu}$ by $(\tilde{X}_t)$.
Let us also denote the Glauber dynamics started at a particular state $x$ by $(X_t^x)$
and $(\tilde{X}_t^x)$ respectively.
The FKG condition implies that if $x \leq y$ then under the monotone coupling (which
will be described in Section \ref{sec:general_tilted}) we also have $X_1^x \leq \tilde{X}_1^y$ almost surely.
Thus we may start both the standard and the tilted chains at the same starting location, say $x$,
and we find that at all times $t$ we have $X_t^x \leq \tilde{X}_t^x$; then just running both
chains to stationarity leads to the stochastic domination result.
Note that the above application of the FKG condition only requires us to consider
states near $x$ and $y$, i.e.\ reachable within one step of the dynamics; thus, the weaker
local FKG condition actually already suffices for the proof of the standard FKG inequality.

In our setting, when conditioning on a metastable well,
we cannot simply run both chains to stationarity, as the monotone coupling may not preserve
ordering when the chains are near the boundary of the well.
This is why we need the mixing condition, so that we may assert that the distributions
of $X_t^x$ and $\tilde{X}_t^x$ are close to $\mu$ and $\tilde{\mu}$ respectively, after
a relatively short amount of time, before either chain has come close to the boundary of
the well.
In fact, to translate the mixing of the untilted dynamics to that for the tilted dynamics,
we require the slightly stronger notion of contraction, which implies rapid mixing.

One important model which is not covered by our current work is the lattice
Ising model which has recently been proven to exhibit a form of rapid metastable
mixing at low temperatures within a phase \cite{gheissari2022low}.
The reason we do not pursue this presently is that the Glauber dynamics for the lattice Ising
model does not exhibit one-step contraction, and the mechanism for mixing is somewhat
more subtle.
Nevertheless, perhaps by considering a block-spin update process, we expect
our methods to be useful in that setting and we hope to pursue this in future work.
\subsection{Outline of the article}
\label{sec:intro_outline}

In Section \ref{sec:general}, we present the general Theorem \ref{thm:general} which formalizes
the ideas of the previous section into a result which may be applied in a variety of situations
to obtain almost nonnegative correlations of increasing functions under the local FKG condition
plus contraction.

In Section \ref{sec:cw}, we apply Theorem \ref{thm:general} to the generalized Curie--Weiss
models introduced in Section \ref{sec:intro_setup_gcwm} and derive Theorem \ref{thm:main_gcwm}.
As we are not aware of any treatment in the literature of GCWMs in the generality which we consider,
this section spells out all of the details and should be a helpful reference
for those wishing to apply Theorem \ref{thm:general} in other settings.

In Section \ref{sec:ergm}, we apply Theorem \ref{thm:general} to the exponential random
graph models introduced in Section \ref{sec:intro_setup_ergm} and derive Theorem
\ref{thm:main_ergm}.
For this application, the inputs to Theorem \ref{thm:general} already exist in the literature
and so we simply recount them in order to finish the proof.

Finally, in Section \ref{sec:covariance}, we present an extension of a classical covariance
inequality using the approximate FKG inequality as an input, and use
this to prove Theorem \ref{thm:main_clt}.
We also mention a variety of other central limit theorems for ERGMs which follow in much the
same way via Theorem \ref{thm:main_ergm}.
\subsection{Acknowledgements}
\label{sec:intro_acknowledgements}

SM would like to thank Subhroshekhar Ghosh and Somabha Mukherjee  for helpful comments on a draft of this paper.
VW would like to thank Subhroshekhar Ghosh and Matteo Mucciconi for an invitation to visit
the National University of Singapore in the summer of 2025, which is where this collaboration
started.
\section{A general theorem for deriving approximate FKG inequalities}
\label{sec:general}

We now state our general result giving almost nonnegative correlations for coordinate-wise increasing functions.
For improved applicability, we have condensed the hypotheses to a few statements which should be checkable in a diverse range of applications.
Moreover, although in the sequel we only consider spin systems on $\{0,1\}^N$, we state our theorem for more general state spaces $\calA^N$,
where $\calA$ is any finite totally ordered set.
While we do not expect finiteness of the spin space to be a strict requirement, our proof does use it in a nontrivial way and it is an interesting
question to get around this for certain applications such as the $\phi^4$ model or other continuum-spin systems.

In what follows, we use $\dh$ to denote the \emph{Hamming distance} on $\calA^N$, which is defined by
\begin{equation}
    \dh(x,y) = \sum_{i=1}^N \ind{x(i) \neq y(i)},
\end{equation}
in other words it is the number of coordinates at which $x$ and $y$ differ.
In addition, for a set $\Lambda \subseteq \calA^N$ and a point $z \in \calA^N$, $\dh(z,\Lambda)$ denotes the minimum distance from $z$ to 
any point in $\Lambda$.

We also introduce the notion of the \emph{intrinsic diameter} $\diam{\Lambda}$ of a set $\Lambda \subseteq \calA^N$.
To define this, let the \emph{intrinsic distance} $\dl(a,b)$ be the length of the shortest
path between $a$ and $b$ \emph{which is contained entirely within $\Lambda$} and takes steps of Hamming length $1$.
Then the intrinsic diameter $\diam{\Lambda}$ is the maximum intrinsic distance between any two points in $\Lambda$.
Note that if $\Lambda$ is not connected in the nearest-neighbor Hamming graph, then $\diam{\Lambda} = \infty$.
It is therefore crucial that to apply the following theorem, the set $\Lambda$ must be connected in this sense.
The key reason for this is the use of a path coupling argument in the proof.

Finally, we let $\Lambda^\comp$ denote the complement $\calA^N \setminus \Lambda$,
and we use $\| f \|_\infty$ to denote the maximum possible absolute value of a function $f : \calA^N \to \R$.
We may now state our general result.

\begin{theorem}
\label{thm:general}
Suppose there is a set $\Lambda \subseteq \calA^N$ which satisfies the following two properties:
\begin{enumerate}[label=(\roman*), leftmargin=*]
    \item
    \label{cond:localFKG}
    The measure satisfies a local version of the FKG condition near $\Lambda$: for any $x,y \in \calA^N$ with both $\dh(x,\Lambda) \leq 1$ and
    $\dh(y,\Lambda) \leq 1$, if we also have $\dh(x \vee y, \{x,y\}) \leq 1$ and $\dh(x \wedge y, \{x,y\}) \leq 1$, then
    \begin{equation}
        \mu(x \vee y) \mu(x \wedge y) \geq \mu(x) \mu(y).
    \end{equation}

    \item
    \label{cond:contraction}
    The Glauber dynamics exhibits one-step contraction under the monotone coupling (to be defined in Section \ref{sec:general_tilted} below)
    when started at neighboring states within $\Lambda$:
    there is some $\alpha < 1$ such that if $a,b \in \Lambda$ and $\dh(a,b) = 1$, then
    \begin{equation}
        \E \left[ \dh(X_1^a, X_1^b) \right] \leq \alpha.
    \end{equation}
\end{enumerate}
Then for any coordinate-wise increasing functions $f,g : \calA^N \to \R$ and any positive integer $T$, we have
\begin{equation}\label{eq:general}
    \Cov_\mu[f,g] \geq - \| f \|_\infty \| g \|_\infty \cdot
    300 \sqrt{T} \left( T \cdot \mu(\Lambda^\comp) + \left( \alpha + \frac{10|\calA|}{T} \right)^T \cdot \diam{\Lambda} \right).
\end{equation}
\end{theorem}

In many situations including our applications later in the article, $\mu(\Lambda^\comp)$
will typically be bounded by $\Exp{-\poly(N)}$ and $\alpha$ will typically have the form $\Exp{- \frac{1}{\poly(N)}}$.
In addition, the intrinsic diameter of $\Lambda$ will also be only $\poly(N)$, and so by taking $T$ to be a particular polynomial (chosen to overcome the
inverse polynomial in the exponent of $\alpha$), we will obtain some form of possibly stretched exponential decay in the error factor.
Since the functions $f$ and $g$ will also typically be only polynomially large, this will suffice for a wide range of applications.

Next, in Section \ref{sec:general_tilted} we introduce some notation and prove Theorem \ref{thm:general} under the assumption of some key bounds,
which will themselves be proved in Section \ref{sec:general_proof} below.

\subsection{Tilted measure and dynamics}
\label{sec:general_tilted}

As mentioned in Section \ref{sec:intro_iop}, our strategy will be to consider a \emph{tilted} version of the measure $\mu$,
as well as the Glauber dynamics with respect to this tilt.
The tilting factor is defined in terms of one of the two functions $f$ and $g$ for which we would like to prove the covariance
lower bound.
However, unlike in Holley's proof of the standard FKG inequality, we require some control on the mixing of both chains, and so
we will require the tilt be mild enough that the mixing properties of the untilted dynamics can be carried over to the tilted setting.

To do this, let us fix some $\eps \in (0,1)$ to be chosen later (it will end up being $\frac{1}{T}$) and define a new function
\begin{equation}
    \tilde{g}(x) \coloneqq g(x) + \frac{2}{\sqrt{\eps}} \| g \|_\infty,
\end{equation}
which satisfies $\tilde{g}(x) > 0$ for all $x \in \calA^N$ since $\eps < 1$,
and also satisfies
\begin{equation}
\label{eq:mildness}
    e^{- \eps} \leq \frac{\tilde{g}(x)}{\E_\mu[\tilde{g}]} \leq e^\eps.
\end{equation}
To see why \eqref{eq:mildness} holds, note that since $e^{s-s^2} \leq 1 + s \leq e^s$ for $|s| < \frac{1}{2}$,
we have
\begin{align}
    \frac{\tilde{g}(x)}{\E_\mu[\tilde{g}]}
    \leq \frac{\left(1 + \frac{2}{\sqrt{\eps}}\right) \| g \|_\infty}{\left(-1 + \frac{2}{\sqrt{\eps}} \right) \| g \|_\infty}
    = \frac{\left(1 + \frac{\sqrt{\eps}}{2}\right)}{\left(1 - \frac{\sqrt{\eps}}{2} \right)}
    \leq \Exp{\frac{\sqrt{\eps}}{2} - \frac{\sqrt{\eps}}{2} + \frac{\eps}{4}},
\end{align}
and similarly for the other inequality in \eqref{eq:mildness}.
Now let us define the tilted measure on $\calA^N$ as
\begin{equation}
    \tilde{\mu}(x) \coloneqq \frac{\tilde{g}(x)}{\E_\mu[\tilde{g}]} \mu(x).
\end{equation}
Then, since $\tilde{g}$ differs from $g$ only by a constant shift, we have
\begin{equation}
\label{eq:covtilt}
    \Cov_\mu[f,g] = \Cov_\mu[f, \tilde{g}] = \E_\mu[\tilde{g}] \left( \E_{\tilde{\mu}}[f] - \E_\mu[f] \right).
\end{equation}
Thus we will obtain almost nonnegative covariance between $f$ and $g$ if we can show that $X \sim \mu$ and $\tilde{X} \sim \tilde{\mu}$ may be coupled so that $f(\tilde{X}) - f(X)$
is almost always positive.
Since $f$ is increasing, for this it suffices to show approximate stochastic dominance of $\mu$ by $\tilde{\mu}$, i.e.\ give a
coupling between $\tilde{X} \sim \tilde{\mu}$ and $X \sim \mu$ which satisfies $\tilde{X} \geq X$ with high probability.

As discussed in Section \ref{sec:intro_iop}, we will achieve this through the \emph{monotone coupling}, which is a simultaneous
coupling of the following four Markov chains:
\begin{align}
    (X_t) &\qquad \text{the stationary $\mu$-Glauber dynamics}, \\
    (\tilde{X}_t) & \qquad \text{the stationary $\tilde{\mu}$-Glauber dynamics}, \\
    (X_t^x) &\qquad \text{the $\mu$-Glauber dynamics started at $X_0^x = x$}, \\
    (\tilde{X}_t^x) &\qquad \text{the $\tilde{\mu}$-Glauber dynamics started at $\tilde{X}_0^x = x$}.
\end{align}
At each time step, the same uniformly random coordinate will be updated in all four chains.
Given that a coordinate $i$ is chosen to be updated, the stationary $\mu$-Glauber dynamics must choose a new value $\sigma \in \calA$
for $X_{t+1}(i)$, keeping the rest of the coordinates the same, with the correct conditional distribution
\begin{equation}
    \mu\left(X(i) = \sigma \middle| X(j) = X_t(j) \text{ for } j \neq i \right),
\end{equation}
and similarly for the other three chains.
In other words, each chain must sample an element of $\calA$ from a particular distribution at each step.
We couple these four distributions at each step under the monotone coupling of measures on the ordered finite set $\calA$,
for instance by sampling a single uniform random number $U$ between $0$ and $1$ and setting the values of the samples from each
distribution according to their inverse CDF at $U$.
This defines the monotone coupling of the chains.

The strategy is now as follows.
We will start $(X_t^z)$ and $(\tilde{X}_t^z)$ at the same state $z$, and start $(X_t)$ and $(\tilde{X}_t)$ from their stationary distributions,
$\mu$ and $\tilde{\mu}$ respectively.
We will then run the chain for time $T$, by which point ideally $(X_t^z)$ and $(X_t)$ will meet so that $X_T^z = X_T$,
and similarly $\tilde{X}_T^z = \tilde{X}_T$; we will use contraction here to show that the monotone coupling forces both non-stationary chains
to meet their corresponding stationary chain.
Simultaneously, the monotone coupling should preserve the order $X_t^z \leq \tilde{X}_t^z$, which will use the local FKG condition
and the fact that $g$ is increasing.
Thus we will obtain $X_T \leq \tilde{X}_T$, and the stationary chains still have distributions $\mu$ and $\tilde{\mu}$ at time $T$,
giving the stochastic domination.
This strategy could run into a problem if any of the chains attempt to leave the set $\Lambda$, but we will show that this is unlikely
as long as we choose the starting point $z$ correctly, and we obtain the following proposition.

\begin{proposition}
\label{prop:keybounds}
Suppose the hypotheses of Theorem \ref{thm:general} hold, as well as the mildness condition \ref{eq:mildness}.
For all $T < \frac{1}{\mu(\Lambda^\comp)}$, there is some $z \in \Lambda$ such that
the following bounds hold under the monotone coupling:
\begin{align}
\label{eq:basechainmixes}
    \P \left[ X_T^z \neq X_T \right] &\leq 2 T \cdot \mu(\Lambda^\comp) + \alpha^T \cdot \diam{\Lambda}, \\
\label{eq:tiltedchainmixes}
    \P \left[ \tilde{X}_T^z \neq \tilde{X}_T \right] &\leq 2 e^{2 \eps T} T \cdot \mu(\Lambda^\comp) + \left( \alpha + 10 |\calA| \eps \right)^T \cdot \diam{\Lambda}, \\
\label{eq:domination}
    \P \left[ \tilde{X}_T^z \not \geq X_T^z \right] &\leq 2 e^{2 \eps T} T \cdot \mu(\Lambda^\comp).
\end{align}
\end{proposition}

We will prove Proposition \ref{prop:keybounds} in Section \ref{sec:general_proof} below,
but before that let us see how we can use the three key bounds it presents to finish the proof of Theorem \ref{thm:general}.

\begin{proof}[Proof of Theorem \ref{thm:general}]
If $T \geq \frac{1}{\mu(\Lambda^\comp)}$ then the right-hand side in the theorem is less than $-300 \| f \|_\infty \| g \|_\infty$,
so the bound holds trivially.
Thus we may assume that $T < \frac{1}{\mu(\Lambda^\comp)}$.
Let $z \in \Lambda$ be the starting point given by Proposition \ref{prop:keybounds}.
Notice that, since the endpoint $\tilde{X}_T$ of the stationary tilted chain has the distribution $\tilde{\mu}$, we have
\begin{equation}
    \E_{\tilde{\mu}}[f] = \E \left[ f(\tilde{X}_T) \right] = \E \left[ f(\tilde{X}_T^z) + \left( f ( \tilde{X}_T ) - f(\tilde{X}_t^z) \right) \ind{\tilde{X}_T^z \neq \tilde{X}_T} \right].
\end{equation}
So, bounding $f(\tilde{X}_T) - f(\tilde{X}_T^z)$ by $2 \| f \|_\infty$, we find that    
\begin{equation}
\label{eq:genbound1}
    \E \left[ f(\tilde{X}_T) \right] \geq \E \left[ f(\tilde{X}_T^z) \right] - 2 \| f \|_\infty \P\left[ \tilde{X}_T^z \neq \tilde{X}_T \right].
\end{equation}
We can do a similar manipulation to replace $\tilde{X}_T^z$ by $X_T^z$, but this time with an inequality; noting that
\begin{equation}
    f(\tilde{X}_T^z) \geq f(X_T^z) + \left( f(\tilde{X}_T^z) - f(X_T^z) \right) \ind{\tilde{X}_T^z \not \geq X_T^z}
\end{equation}
by the increasing nature of $f$, we obtain
\begin{equation}
\label{eq:genbound2}
    \E \left[ f(\tilde{X}_T^z) \right] \geq \E \left[ f(X_T^z) \right] - 2 \|f\|_\infty \P \left[ \tilde{X}_T^z \not \geq X_T^z \right].
\end{equation}
Finally, by the same manipulation as \eqref{eq:genbound1}, we find that
\begin{equation}
\label{eq:genbound3}
    \E \left[ f(X_T^z) \right] \geq \E \left[ f(X_T) \right] - 2 \| f \|_\infty \P \left[ X_T^z \neq X_T \right].
\end{equation}
Since $X_T$ has the distribution $\mu$, by combining \eqref{eq:genbound1}, \eqref{eq:genbound2}, and \eqref{eq:genbound3} we obtain
\begin{equation}
    \E_{\tilde{\mu}}[f] - \E_\mu[f] \geq - 2 \| f \|_\infty
    \left( \P \left[ \tilde{X}_T^z \neq \tilde{X}_T \right] + \P \left[ \tilde{X}_T^z \not \geq X_T^z \right] + \P \left[ X_T^z \neq X_T \right] \right).
\end{equation}
Plugging in the bounds from Proposition \ref{prop:keybounds} above (and replacing \eqref{eq:basechainmixes} with \eqref{eq:tiltedchainmixes} as it is larger), we find that
\begin{equation}
    \E_{\tilde{\mu}}[f] - \E_\mu[f] \geq - 2 \| f \|_\infty \left( 6 e^{2 \eps T} T \cdot \mu(\Lambda^\comp) + 2 \left( \alpha + 10 |\calA| \eps \right)^T \cdot \diam{\Lambda} \right).
\end{equation}
Thus, choosing $\eps = \frac{1}{T}$ and using the fact that $e^2 < 8$ we find that
\begin{equation}
    \E_{\tilde{\mu}}[f] - \E_\mu[f] \geq - 2 \| f \|_\infty \left( 48 T \cdot \mu(\Lambda^\comp) + 2 \left( \alpha + \frac{10 |\calA|}{T} \right)^T \cdot \diam{\Lambda} \right).
\end{equation}
Finally, by \eqref{eq:covtilt} and since $0 < \E_\mu[\tilde{g}] \leq \frac{3}{\sqrt{\eps}} \| g \|_\infty = 3\sqrt{T} \| g \|_\infty$, we find that
\begin{equation}
    \Cov_\mu[f,g] \geq - \| f \|_\infty \| g \|_\infty \cdot 300 \sqrt{T} \left( T \cdot \mu(\Lambda^\comp) + \left( \alpha + \frac{10|\calA|}{T} \right)^T \cdot \diam{\Lambda} \right).
\end{equation}
This finishes the proof.
\end{proof}
\subsection{Proofs of key bounds}
\label{sec:general_proof}

In this section we prove Proposition \ref{prop:keybounds}, finishing the proof of Theorem \ref{thm:general}.
Throughout this section we assume the hypotheses of Theorem \ref{thm:general} hold.
Moreover, we assume that the tilt mildness condition \eqref{eq:mildness} holds with parameter $\eps$.
As a first step in our derivation of Proposition \ref{prop:keybounds}, we remark that,
as a trivial consequence of the mildness condition, the complement of $\Lambda$ is still small under the tilted measure.

\begin{lemma}
\label{lem:tilted_lambdalarge}
$\tilde{\mu}(\Lambda^\comp) \leq e^{\eps} \cdot \mu(\Lambda^\comp)$.
\end{lemma}

Next, we introduce the starting point $z \in \Lambda$ which will be used in Proposition \ref{prop:keybounds}.
It turns out that all we need, given the hypotheses in Theorem \ref{thm:general}, is for both of the chains started at $z$ to remain in $\Lambda$
for enough time, with high probability.
The existence of $z \in \Lambda$ satisfying this condition follows from the condition that $T < \frac{1}{\mu(\Lambda^\comp)}$ by a union bound,
which will be carried out in the proof of the next lemma.

\begin{lemma}
\label{lem:starting_point}
For all $T < \frac{1}{\mu(\Lambda^\comp)}$, there is some $z \in \Lambda$ such that
\begin{align}
    \P \left[ X_t^z \in \Lambda \text{ for } 0 \leq t < T \right] &\geq 1 - T \cdot \mu(\Lambda^\comp), \\
    \P \left[ \tilde{X}_t^z \in \Lambda \text{ for } 0 \leq t < T \right] &\geq 1 - e^{2 \eps T} T \cdot \mu(\Lambda^\comp).
\end{align}
\end{lemma}

\begin{proof}[Proof of Lemma \ref{lem:starting_point}]
Let us begin with the first inequality.
Since $(X_t)$ is stationary with distribution $\mu$, by a union bound we have
\begin{equation}
    \P \left[ X_t \notin \Lambda \text{ for some } 0 \leq t < T \right] \leq T \cdot \mu(\Lambda^\comp).
\end{equation}
Conditioning on $X_0$, we find that
\begin{equation}
    \E \left[ \P \left[ X_t \notin \Lambda \text{ for some } 0 \leq t < T \middle| X_0 \right] \right] \leq T \cdot \mu(\Lambda^\comp).
\end{equation}
Thus the inner probability above must take a value $\leq T \cdot \mu(\Lambda^\comp)$ for some choice of $X_0 = z$,
and we choose this $z$ as our starting point.
Note that we must have $z \in \Lambda$ since otherwise the inner probability is $1$ and we have assumed that $T \cdot \mu(\Lambda^\comp) < 1$.

Now for the second inequality, we cannot simply repeat the above strategy for the tilted chain since we may not arrive at the
\emph{same} starting point $z$ for both chains.
Of course, there are ways to remedy this, but instead we choose to use the mildness assumption \eqref{eq:mildness} directly,
as it will lead us to a bound which will also be useful in the sequel.

We will show below that the update probabilities differ in the tilted dynamics (as compared to the base dynamics) by a factor of at most $e^{2\eps}$ at each step.
Thus any trajectory of length $T$ has probability at most $e^{2 \eps T}$ times larger under the tilted dynamics
as compared to the base dynamics.
Consider the set of all trajectories starting at $z$ and which eventually leave $\Lambda$ by time $T$.
Since the set of these trajectories has probability at most $T \cdot \mu(\Lambda^\comp)$ under the base dynamics, it has probability at most $e^{2 \eps T} \cdot T \cdot \mu(\Lambda^\comp)$
under the tilted dynamics, finishing the proof.

It just remains to prove that the update probabilities differ by at most $e^{2\eps}$ at each step.
Suppose that we start at the state $a$ and that the coordinate $i$ is chosen to be updated; let us consider the probability of arriving at $a^{\sigma \to i}$, which is
our notation for the state $a$ but with $a(i)$ set to $\sigma$ instead of whatever it was before.
By the definition of Glauber dynamics for both the base and tilted chains, as well as the mildness condition \eqref{eq:mildness}, we have
\begin{equation}
\label{eq:update_UB}
    \P\left[ \tilde{X}_1^a = a^{\sigma \to i} \middle| i \text{ chosen} \right] = \frac{\tilde{\mu}(a^{\sigma \to i})}{\sum_{\tau \in \calA} \tilde{\mu}(a^{\tau \to i})}
    \leq \frac{e^{\eps} \mu(a^{\sigma \to i})}{\sum_{\tau \in \calA} e^{-\eps} \mu(a^{\tau \to i})}
    = e^{2\eps} \P \left[ X_1^a = a^{\sigma \to i} \middle| i \text{ chosen} \right],
\end{equation}
and similarly we have
\begin{equation}
\label{eq:update_LB}
    \P \left[ \tilde{X}_1^a = a^{\sigma \to i} \middle| i \text{ chosen} \right] \geq e^{-2\eps} \P \left[ X_1^a = a^{\sigma \to i} \middle| i \text{ chosen} \right].
\end{equation}
Thus, averaging over the uniform choice of coordinate $i$ to update in both of the Glauber dynamics finishes the proof.
\end{proof}

One more ingredient is needed for the proof of Proposition \ref{prop:keybounds}: we will use the mildness condition \eqref{eq:mildness},
and in particular the corresponding fact for update probabilities just described in the previous proof, to derive a one-step contraction bound
under the monotone coupling for the tilted chain from the corresponding fact (condition \ref{cond:contraction}) for the base chain.

\begin{lemma}
\label{lem:tilted_contraction}
If $a, b \in \Lambda$ differ in only one entry, then under the monotone coupling we have
\begin{equation}
    \E  \left[ \dh(\tilde{X}_1^a, \tilde{X}_1^b) \right] \leq \alpha + 10 |\calA| \eps.
\end{equation}
\end{lemma}

\begin{proof}[Proof of Lemma \ref{lem:tilted_contraction}]
Suppose that $a$ and $b$ differ at coordinate $j$.
Then under the monotone coupling, if $j$ is chosen to be updated, we will obtain $\tilde{X}_1^a = \tilde{X}_1^b$.
If instead some other coordinate $i \neq j$ is chosen to be updated, then $\dh(\tilde{X}_1^a, \tilde{X}_1^b) \geq 1$, and it may equal $2$
if we subsequently have $\tilde{X}_1^a(i) \neq \tilde{X}_1^b(i)$.
We need to control the probability that this happens.

First let us suppose that we have any two distributions $\nu$ and $\nu'$ on $\calA$; later, we will take these to be the distributions of $\tilde{X}_1^a(i)$ and $\tilde{X}_1^b(i)$
in one case, or the distributions of $X_1^a(i)$ and $X_1^b(i)$ in another, both conditioned on choosing coordinate $i$ to update in the Glauber dynamics.
Let us denote by $(\tau,\tau') \in \calA^2$  the random variable sampled from the monotone coupling between $\nu$ and $\nu'$.
Recall that this is constructed by first sampling $U \in [0,1]$ uniformly, and setting $\tau = \sigma$ if $\nu(\tau < \sigma) \leq U \leq \nu(\tau \leq \sigma)$ and likewise for $\tau'$,
using the same sample $U$.
The probability that $\tau = \tau'$ is thus the probability that $U$ lands in a collection of ``good'' intervals where the inverse CDFs of $\nu$ and $\nu'$ agree.
The endpoints of these intervals are of the form $\nu(\tau \leq \sigma)$ or $\nu'(\tau' \leq \sigma')$ for some $\sigma$ or $\sigma'$ in $\calA$.

Now let us start with $\nu$ being the distribution of $X_1^a(i)$ and $\nu'$ being the distribution of $X_1^b(i)$, given that $i$ is chosen to be updated.
The total length of the good intervals for these distributions is exactly 
\begin{equation}
    \P \left[ X_1^a(i) = X_1^b(i) \middle| i \text{ chosen} \right],
\end{equation}
and the endpoints of the good intervals are of the form
\begin{equation}
    \P \left[ X_1^a(i) \leq \sigma \middle| i \text{ chosen} \right]
    \qquad \text{or} \qquad
    \P \left[ X_1^b(i) \leq \sigma' \middle| i \text{ chosen} \right],
\end{equation}
for some $\sigma$ or $\sigma'$ in $\calA$.

Now if we switch $\nu$ and $\nu'$ to the distributions of $\tilde{X}_1^a(i)$ and $\tilde{X}_1^b(i)$ (again given $i$ is chosen), then according to \eqref{eq:update_UB}
and \eqref{eq:update_LB}, the endpoints of the good intervals can each shift by only a small amount.
Indeed, note first that \eqref{eq:update_UB} easily implies
\begin{equation}
\label{eq:8eps}
    \P \left[ \tilde{X}_1^a(i) \leq \sigma \middle| i \text{ chosen} \right] \leq e^{2\eps} \P \left[ X_1^a(i) \leq \sigma \right]
    \leq \P \left[ X_1^a(i) \leq \sigma \right] + 8 \eps,
\end{equation}
using the inequality $e^{2\eps} \leq 1 + 8 \eps$ for $\eps \in (0,1)$ and the fact that probabilities are at most $1$.
Similarly, \eqref{eq:update_LB} implies that
\begin{equation}
\label{eq:2eps}
    \P \left[ \tilde{X}_1^a(i) \leq \sigma \middle| i \text{ chosen} \right] \geq e^{-2\eps} \P \left[ X_1^a(i) \leq \sigma \right]
    \geq \P \left[ X_1^a(i) \leq \sigma \right] - 2 \eps,
\end{equation}
using the inequality $e^{-2\eps} \geq 1-2\eps$ this time.
Of course, we also have the inequalities \eqref{eq:8eps} and \eqref{eq:2eps} if $a$ is replaced by $b$.
Thus the lower endpoint of a good interval may shift up by at most $8 \eps$ and the upper endpoint of a good interval may shift down by at most $2 \eps$.
So the total length of each good interval can decrease by at most $10 \eps$.
There are at most $|\calA|$ good intervals, so we find that
\begin{equation}
    \P \left[ \tilde{X}_1^a(i) = \tilde{X}_1^b(i) \middle| i \text{ chosen} \right] \geq \P \left[ X_1^a(i) = X_1^b(i) \middle| i \text{ chosen} \right] - 10 |\calA| \eps.
\end{equation}
Averaging over the choice of $i$, we find that
\begin{equation}
    \P \left[ \dh(\tilde{X}_1^a, \tilde{X}_1^b) = 2 \right] \leq \P \left[ \dh(X_1^a, X_1^b) = 2 \right] + 10 |\calA| \eps.
\end{equation}
Now recall that the only values the Hamming distance can take at this stage are $0, 1$ and $2$, and both $\dh(\tilde{X}_1^a,\tilde{X}_1^b)$ and $\dh(X_1^a,X_1^b)$
take the value $0$ with the same probability, which is exactly the probability that $j$ (the coordinate where $a$ and $b$ differ) is chosen to be updated.
This also implies that the probability that each Hamming distance is $\geq 1$ is also the same, so we find that
\begin{align}
    \E \left[ \dh(\tilde{X}_1^a,\tilde{X}_1^b) \right] &= \P \left[ \dh(\tilde{X}_1^a, \tilde{X}_1^b) \geq 1 \right] + \P \left[ \dh(\tilde{X}_1^a,\tilde{X}_1^b) = 2 \right] \\
    &\leq \P \left[ \dh(X_1^a, X_1^b) \geq 1 \right] + \P \left[ \dh(X_1^a,X_1^b) = 2 \right] + 10 |\calA| \eps \\
    &= \E \left[ \dh(X_1^a,X_1^b) \right] + 10 |\calA| \eps.
\end{align}
So by condition \ref{cond:contraction} which states that the last expectation above is $\leq \alpha$, the proof is finished.
\end{proof}

With Lemmas \ref{lem:tilted_lambdalarge}, \ref{lem:starting_point}, and \ref{lem:tilted_contraction} in hand, we are now ready to provide the key bounds used in the
proof of Theorem \ref{thm:general}, which were stated above in Proposition \ref{prop:keybounds}.
Recall that this states that for $T < \frac{1}{\mu(\Lambda^\comp)}$ there is some $z \in \Lambda$ such that the chains started at $z$ mix to stationarity rapidly, and such that the tilted
chain started at $z$ remains above the base chain started at $z$ up to time $T$.

\begin{proof}[Proof of Proposition \ref{prop:keybounds}]
Let $z \in \Lambda$ be the starting point given by Lemma \ref{lem:starting_point} above; note that the existence of this point
is the only place where we use the assumption $T < \frac{1}{\mu(\Lambda^\comp)}$.
We will quickly prove the first two inequalities, \eqref{eq:basechainmixes} and \eqref{eq:tiltedchainmixes}, which prove that the chains mix rapidly by giving upper bounds on
$\P \left[ X_T^z \neq X_T \right]$ and $\P \left[ \tilde{X}_T^z \neq \tilde{X}_T \right]$.
The inequality \eqref{eq:domination}, which gives an upper bound on $\P \left[ \tilde{X}_T^z \not \geq X_T^z \right]$ is a bit trickier, and this is where we will use
the local FKG condition \ref{cond:localFKG} in the hypotheses of Theorem \ref{thm:general}.

Let us focus first on the bound \eqref{eq:basechainmixes} for the standard dynamics.
Recalling the definition of the intrinsic distance $\dl$ from before the statement of Theorem \ref{thm:general}, note that we have
\begin{equation}
    \P \left[ X_T^z \neq X_T \right] \leq \E \left[ \dl(X_T^z, X_T) \right].
\end{equation}
By a union bound as in the proof of Lemma \ref{lem:starting_point}, with probability at least $1 - T \cdot \mu(\Lambda^\comp)$,
the stationary chain remains in $\Lambda$ for $0 \leq t < T$.
Additionally, the non-stationary chain $X_t^z$ remains in $\Lambda$ for $0 \leq t < T$ with probability at least $1 - T \cdot \mu(\Lambda^\comp)$ by Lemma \ref{lem:starting_point}.
On the event where both chains remain in $\Lambda$ for the entire duration up to time $T$, the contraction of condition \ref{cond:contraction} yields reduction
of the expected intrinsic distance between the two chains by a factor of $\alpha$ for each time step, by a standard path coupling argument.
Since $X_0^z$ and $X_0$ are at most $\diam{\Lambda}$ apart in $\Lambda$-distance if they both start in $\Lambda$, we thus obtain
\begin{equation}
    \P \left[ X_T^z \neq X_T \right] \leq 2 T \cdot \mu(\Lambda^\comp) + \alpha^T \cdot \diam{\Lambda}.
\end{equation}
The same argument applies to the tilted dynamics, but we must apply Lemma \ref{lem:tilted_lambdalarge} for the bound on $\tilde{\mu}(\Lambda^\comp)$.
Also, instead of condition \ref{cond:contraction} we must apply Lemma \ref{lem:tilted_contraction}, and we must use the second part of Lemma \ref{lem:starting_point}
rather than the first.
With these modifications, we obtain
\begin{align}
    \P \left[ \tilde{X}_T^z \neq \tilde{X}_T \right] &\leq (e^\eps + e^{2 \eps T}) T \cdot \mu(\Lambda^\comp) + \left( \alpha + 10 |\calA| \eps \right)^T \cdot \diam{\Lambda} \\
    &\leq 2 e^{2 \eps T} T \cdot \mu(\Lambda^\comp) + \left( \alpha + 10 |\calA| \eps \right)^T \cdot \diam{\Lambda}
\end{align}
as desired, proving the bound \eqref{eq:tiltedchainmixes} as well.

Now we turn to the proof of the bound \eqref{eq:domination}, showing that $\tilde{X}_T^z \not \geq X_T^z$ is unlikely.
Suppose that $a, b \in \Lambda$ and $a \leq b$.
Then given that coordinate $i$ is chosen to be updated, for any $\sigma \in \calA$ we have
\begin{equation}
\label{eq:rhsdec}
    \P \left[ \tilde{X}_1^b(i) \geq \sigma \middle| i \text{ chosen} \right]
    = \frac{\sum_{\tau \geq \sigma} \tilde{\mu} (b^{\tau \to i})}{\sum_{\tau \in \calA} \tilde{\mu} (b^{\tau \to i}) } 
    = \frac{\sum_{\tau \geq \sigma} \tilde{g}(b^{\tau \to i}) \mu (b^{\tau \to i})}{\sum_{\tau \in \calA} \tilde{g}(b^{\tau \to i}) \mu (b^{\tau \to i}) }.
\end{equation}
Now, since $\tilde{g}$ is increasing, the right-hand side above cannot increase if we remove the factors of $\tilde{g}$.
To see why this holds in more detail, first note that
if $r \geq s > 0$ and $p, q \geq 0$ with $p + q > 0$ then we have
\begin{equation}
\label{eq:rprpsq}
    \frac{rp}{rp + sq} \geq \frac{p}{p + q}.
\end{equation}
We will take
\begin{equation}
    p \coloneqq \sum_{\tau \geq \sigma} \mu(b^{\tau \to i})
    \qquad \text{and} \qquad
    q \coloneqq \sum_{\tau < \sigma} \mu(b^{\tau \to i}).
\end{equation}
Note that if either one of these is zero, then removing the factors of $\tilde{g}$ from the right-hand side of \eqref{eq:rhsdec} has no effect, i.e.\ it cannot increase.
So let us assume that neither $p$ nor $q$ is zero, and define
\begin{equation}
    r \coloneqq \frac{\sum_{\tau \geq \sigma} \tilde{g}(b^{\tau \to i}) \mu(b^{\tau \to i})}{\sum_{\tau \geq \sigma} \mu(b^{\tau \to i})}
    \geq \tilde{g}(b^{\sigma \to i})
    \geq \frac{\sum_{\tau < \sigma} \tilde{g}(b^{\tau \to i}) \mu(b^{\tau \to i})}{\sum_{\tau < \sigma} \mu(b^{\tau \to i})} \eqcolon s.
\end{equation}
Then applying \eqref{eq:rprpsq}
we see that the right-hand side in \eqref{eq:rhsdec} cannot increase if we erase the factors of $\tilde{g}$.
So we find that
\begin{equation}
\label{eq:dom_tilt}
    \P \left[ \tilde{X}_1^b(i) \geq \sigma \middle| i \text{ chosen} \right] \geq \frac{\sum_{\tau \geq \sigma} \mu(b^{\tau \to i})}{\sum_{\tau \in \calA} \mu(b^{\tau \to i})}.
\end{equation}
By definition, we also have
\begin{equation}
\label{eq:dom_base}
    \P \left[ X_1^a(i) \geq \sigma \middle| i \text{ chosen} \right] = \frac{\sum_{\tau \geq \sigma} \mu(a^{\tau \to i})}{\sum_{\tau \in \calA} \mu(a^{\tau \to i})}.
\end{equation}
We wish to show that the right-hand side of \eqref{eq:dom_tilt} is greater than or equal to the right-hand side of \eqref{eq:dom_base}.
This is equivalent to showing that
\begin{equation}
    \left( \sum_{\tau_+ \geq \sigma} \mu(b^{\tau_+ \to i}) \right)  \left( \sum_{\tau_- \in \calA} \mu (a^{\tau_- \to i}) \right)
    \geq \left( \sum_{\tau_+ \geq \sigma} \mu(a^{\tau_+ \to i}) \right)  \left( \sum_{\tau_- \in \calA} \mu (b^{\tau_- \to i}) \right).
\end{equation}
Expanding both sides and erasing the common terms (where $\tau_- \geq \calA$), to prove the above it suffices to show that
\begin{equation}
\label{eq:fkg_need}
    \mu(b^{\tau_+ \to i}) \mu(a^{\tau_- \to i}) \geq \mu(a^{\tau_+ \to i}) \mu(b^{\tau_- \to i}).
\end{equation}
for $\tau_+ \geq \sigma > \tau_-$.
Setting $x = a^{\tau_+ \to i}$ and $y = b^{\tau_- \to i}$, since $a \leq b$ we find that
\begin{equation}
    x \vee y = b^{\tau_+ \to i} \qquad \text{and} \qquad x \wedge y = a^{\tau_- \to i}.
\end{equation}
Thus \eqref{eq:fkg_need} follows by an application of the local version of the FKG condition \ref{cond:localFKG},
since $a$ and $b$ are in $\Lambda$ and $x$ and $y$ are at Hamming distance $1$ from these.
We thus find that if $a,b \in \Lambda$ and $a \leq b$ then
\begin{equation}
    \P \left[ \tilde{X}_1^b(i) \geq \sigma \right] \geq \P \left[ X_1^a(i) \geq \sigma \right]
\end{equation}
for all $\sigma \in \calA$, given that the coordinate $i$ is chosen to be updated.

This implies that if $\tilde{X}_{t-1}^z \geq X_{t-1}^z$ and both of these are in $\Lambda$, then we will also have $\tilde{X}_t^z \geq X_t^z$
under the monotone coupling.
And since both chains start at the same point $z \in \Lambda$, we do have $\tilde{X}_0^z \geq X_0^z$.
Thus we find that
\begin{equation}
    \P \left[ \tilde{X}_T^z \not \geq X_T^z \right] \leq \P \left[ \tilde{X}_t^z \notin \Lambda \text{ or } X_t^z \notin \Lambda \text{ for some } t \leq T \right].
\end{equation}
The right-hand side above is bounded by
\begin{equation}
    \left( 1 + e^{2 \eps T} \right) T \cdot \mu(\Lambda^\comp) \leq 2 e^{2 \eps T} T \cdot \mu(\Lambda^\comp)
\end{equation}
by Lemma \ref{lem:starting_point}.
This proves \eqref{eq:domination}, finishing the proof of Proposition \ref{prop:keybounds} and thus of Theorem \ref{thm:general}.
\end{proof}
\section{Warm-up application to generalized Curie--Weiss models}
\label{sec:cw}

In this section we will apply our general result, Theorem \ref{thm:general} to the generalized Curie--Weiss models introduced
in Section \ref{sec:intro_setup_gcwm}, resulting in a proof of Theorem \ref{thm:main_gcwm}.
Throughout this entire section, for ease of notation, we set $h = h_\beta^\GCW$ to be the polynomial defined in \eqref{eq:polyGCW},
and set $\hamilton = \hamiltonGCW$ to be the Hamiltonian defined in \eqref{eq:hamiltonGCW}, so that $\hamilton(x) = h(m(x))$
for all $x \in \{0,1\}^N$, where the mean magnetization $m(x) = \frac{1}{N} \sum_{i=1}^N x(i)$ was initially defined in \eqref{eq:magnetization}.
We also set $\nu = \nu^\GCW_\beta$ to be the full (unconditioned) measure on $\{0,1\}^N$ proportional to $e^{N \hamilton(x)} = e^{N h(m(x))}$
as defined in \eqref{eq:cw_full}.
Additionally, we set $L = \rateGCW$ as defined in \eqref{eq:rateGCW}, so that
\begin{equation}
\label{eq:Lagain}
    L(m) = h(m) - m \log m - (1-m) \log(1-m)
\end{equation}
for $m \in [0,1]$.
This is the rate function for the large deviations principle which we will prove shortly.

To obtain the contraction estimate used in condition \ref{cond:contraction} of Theorem \ref{thm:general}, it will also be useful
to note an equivalence between convex local maximizers of $L$ and \emph{attracting fixed points} of a function on $[0,1]$ which
controls the dynamics, as will be made clear in Section \ref{sec:cw_contraction} below.
To state this, let us define $\varphi(s) = \frac{e^s}{1+e^s}$.

\begin{lemma}
\label{lem:equiv}
The following two statements are equivalent for any $m \in [0,1]$:
\begin{enumerate}
    \item $L'(m) = 0$ and $L''(m) < 0$,
    \item $m = \varphi(h'(m))$ and $\frac{d}{d m} \varphi(h'(m)) < 1$.
\end{enumerate}
\end{lemma}

\begin{proof}[Proof of Lemma \ref{lem:equiv}]
We have
\begin{equation}
\label{eq:inter}
    L'(m) = h'(m) - \log \frac{m}{1-m},
\end{equation}
so since $\varphi(s) = \frac{e^s}{1+e^s}$ is the inverse of $m \mapsto \log \frac{m}{1-m}$, we have $L'(m) = 0$ exactly when
$m = \varphi(h'(m))$.
Note also that $\varphi'(s) = \varphi(s) (1-\varphi(s))$, which we will use shortly.
Taking another derivative of \eqref{eq:inter} we find that
\begin{equation}
    L''(m) = h''(m) - \frac{1}{m(1-m)},
\end{equation}
so if $L'(m) = 0$ and $L''(m) < 0$ then multiplying the right-hand side above by $\varphi'(h'(m)) = \varphi(h'(m)) (1-\varphi(h'(m))) = m(1-m)$ we find that
\begin{equation}
    \varphi'(h'(m)) h''(m) < 1,
\end{equation}
and the left-hand side is $\frac{d}{dm} \varphi(h'(m))$.
These steps can be reversed, showing the equivalence.
\end{proof}

Theorem \ref{thm:main_gcwm} concerns not the full measure $\nu$ but rather the phase measure $\calM_{\beta,m_*}^\GCW$, which is the measure $\nu$
\emph{conditioned} on the set
\begin{equation}
    \Lambda_{m_*}^\eta \coloneqq \{ x : |m(x) - m_*| \leq \eta \}
\end{equation}
for some fixed $m_* \in U_\beta^\GCW$.
Recall that $U_\beta^\GCW$ denotes the set of global maximizers of $L$ which also satisfy $L''(m_*) < 0$, as introduced in Section
\ref{sec:intro_setup_gcwm}, so that they are also attracting fixed points of $m \mapsto \varphi(h'(m))$ by Lemma \ref{lem:equiv}.
Here $\eta$ is chosen to be small enough that the sets $\Lambda_{m_*}^\eta$ are disjoint as $m_*$ varies in $M_\beta^\GCW$, the
set of \emph{all} global maximizers of $L$.
We will typically suppress $\eta$ from our notation as the results do not depend on $\eta$ once it is small enough.
However, in the course of the proof, we will also make use of $\Lambda_{m_*}^\eps$ for values of $\eps$ smaller than $\eta$,
which motivates the introduction of these sets explicitly.
Finally, for notational convenience, we set $\mu = \calM_{\beta,m_*}^\GCW$, as we will only consider one phase measure at a time.

\subsection{Large deviations principle}
\label{sec:cw_ldp}

In this section we prove the large deviations principle previously mentioned in \eqref{eq:cw_ldp_initial}, validating the decomposition
of the measure $\nu$ into the phase measures $\calM_{\beta,m_*}^\GCW$.
The proofs in this section are standard but we include them for completeness, as we are not aware of a treatment of GCWMs in the literature
as we have defined them, although certain special cases have appeared \cite{liggett2007statistical,samanta2024mixing}.

First, we need the following lemma which shows that the rate function $L$ indeed captures the exponential rate of the probability
under $\nu$ of observing a particular magnetization.

\begin{lemma}
\label{lem:freeenergy}
For any $m \in [0,1]$ which is an integer multiple of $\frac{1}{N}$, we have
\begin{equation}
    \binom{N}{m N} e^{N h(m)} = \Exp{N \cdot L(m) + O(\log N)},
\end{equation}
where the constant in the $O(\log N)$ error does not depend on $m$.
\end{lemma}

\begin{proof}[Proof of Lemma \ref{lem:freeenergy}]
This follows directly from Stirling's approximation in the form
\begin{equation}
    \log(n!) = n \log n - n + O(\log n)
\end{equation}
as $n \to \infty$, plus the special case $\log(0!) = 0 = 0 \log 0 - 0$, with the interpretation that $0 \log 0 = 0$.
Applying this to the logarithm of the expression we wish to approximate, we find that
\begin{align}
    \log \left(
        \binom{N}{m N} e^{N h(m)}
    \right)
    &= N h(m) + \log(N!) - \log((mN)!) - \log(((1-m)N)!) \\
    &= N h(m) + N \log N - N - m N \log (m N) +  m N \\
    &\qquad - (1-m) N \log ((1-m)N) + (1-m)N + O(\log N)\\
    &=  N \left( h(m) - m \log m - (1-m) \log (1-m) \right) + O(\log N)
\end{align}
as required, recalling the formula \eqref{eq:Lagain} for $L$.
\end{proof}

Now we are ready to prove the large deviations principle itself, validating \eqref{eq:cw_ldp_initial}.

\begin{proposition}
\label{prop:gcw_ldp}
For any $\eta > 0$ there are some $C(\eta), c(\eta) > 0$ such that we have
\begin{equation}
    \nu \left( \{ 0, 1 \}^N \setminus \bigcup_{m_* \in M_\beta^\GCW} \Lambda_{m_*}^\eta \right) \leq C(\eta) e^{- c(\eta) N}.
\end{equation}
\end{proposition}

\begin{proof}[Proof of Proposition \ref{prop:gcw_ldp}]
First note that for all $x \in \{0,1\}^N$, $m(x)$ is an integer multiple of $\frac{1}{N}$.
Let us denote by $B$ the set of such integer multiples in $[0,1]$ which are at least $\eta$ away from any maximizer $m_* \in M_\beta^\GCW$.
By the continuity of $L$, there is some $\delta>0$ such that $L(\tilde{m}) < L(m_*) - \delta$ for all $\tilde{m} \in B$.
Thus the probability in the proposition which we would like to bound is at most
\begin{equation}
    N \cdot \max \left\{  \nu \left( \left\{ x : m(x) = \tilde{m} \right\} \right) : \tilde{m} \in B \right\}
    \leq \frac{N}{\calZ} \cdot \max \left\{ \binom{N}{\tilde{m} N} e^{N h(\tilde{m})}  : \tilde{m} \in B \right\},
\end{equation}
where we write $\calZ$ for the \emph{partition function} of $\nu$, i.e.
\begin{equation}
    \calZ \coloneqq \sum_{x \in \{0,1\}^N} e^{N h(m(x))} \geq \binom{N}{m N} e^{N h(m)}
\end{equation}
for any individual $m$ which is an integer multiple of $\frac{1}{N}$.
Let us choose $\hat{m}$ to be the integer multiple of $\frac{1}{N}$
which is as close as possible to some $m_* \in M_\beta^\GCW$.
This means that $|\hat{m} - m_*| < \frac{1}{N}$, so by differentiability of $L$ at $m_*$, we also have
$L(\hat{m}) > L(m_*) - O\left(\frac{1}{N}\right)$.
Note that here we have used the fact that $m_* \notin \{0,1\}$, since $L'(m) \to \infty$ as $m \to 0$ and $L'(m) \to - \infty$ as $m \to 1$.
Thus, applying Lemma \ref{lem:freeenergy}, the probability in the proposition which we would like to bound is at most
\begin{equation}
    N \Exp{ N \left(L(m_*) - \delta\right) - N \left(L(m_*) - O\left(\frac{1}{N}\right)\right) + O(\log N)} 
    = \Exp{- \delta N + O(\log N)},
\end{equation}
finishing the proof.
\end{proof}

Note that the same proof also shows that, once we fix $\eta > 0$ so that the metastable well measures $\mu = \calM_{\beta,m_*}^\GCW$ are well-defined
(owing to the disjointness of $\Lambda_{m_*}^\eta$ for $m_* \in M_\beta^\GCW$),
then for any $\eps \in (0,\eta)$ we have
\begin{equation}
\label{eq:lambdalarge_gcw}
    \mu(\Lambda_{m_*}^\eps) > 1 - e^{-\Omega(N)},
\end{equation}
where the constant in the $\Omega(N)$ term depends on $\eps$.

\subsection{Contraction in almost the entire metastable well}
\label{sec:cw_contraction}

An important ingredient in the application of Theorem \ref{thm:general} is the contraction estimate within some large enough set.
In the present case of GCWMs, this large set will simply be $\Lambda_{m_*}^\eps$ for some $\eps$ smaller than
the value $\eta$ in the definition of the phase measure $\mu = \calM_{\beta,m_*}^\GCW$; note that this will not be the case for ERGMs
considered in Section \ref{sec:ergm} below.
Proposition \ref{prop:gcw_contraction} below establishes the required contraction for GCWMs.

First we remark that the update probabilities for the Glauber dynamics of GCWMs may be written explicitly.
To do this, let us introduce the notion of \emph{discrete derivatives}: for each $i \in \{1,\dotsc,N\}$ and any function
$f : \{0,1\}^N \to \R$ we write $\partial_i f(x) = f(x^{+i}) - f(x^{-i})$, where $x^{+ i}$ denotes the spin configuration $x$
with the $i$th coordinate set to $1$, and $x^{- i}$ denotes the configuration with the $i$th coordinate set to $0$.
Recall that we denote by $(X_t^x)$ the Glauber dynamics started at $X_0^x = x$.
By the definition of the full GCWM measure $\nu$, the probability of ending up at $x^{+i}$ after one step, given that we chose
coordinate $i$ to update, is exactly
\begin{equation}
    \frac{e^{N \hamilton(x^{+i})}}{e^{N \hamilton(x^{-i})} + e^{N \hamilton(x^{+i})}}
    = \varphi( N \partial_i \hamilton(x) ),
\end{equation}
recalling the definition of $\varphi(s) = \frac{e^s}{1+e^s}$ introduced above.
If we start at a state $x \in \Lambda_{m_*}^\eps$ for some $\eps$ \emph{strictly smaller} than the $\eta$ in the definition of the phase measure,
then the \emph{conditioned} Glauber dynamics also obeys the same formula since it will have no chance to escape the larger
set $\Lambda_{m_*}^\eta$.
In any case, to control this quantity, we first present a short lemma relating the discrete derivative of polynomials of $m(x)$,
such as $\hamilton(x) = h(m(x))$, with the standard notion of derivatives of polynomials.

\begin{lemma}
\label{lem:gcw_partial}
For any $i,j \in \{1,\dotsc,N\}$ and any $x \in \{0,1\}^N$, we have
\begin{align}
    \partial_i (m(x))^k &= \frac{k m(x^{-i})^{k-1}}{N} + O \left(\frac{1}{N^2} \right), \\
    \partial_j \partial_i (m(x))^k &= \frac{k(k-1) m (x^{-i-j})^{k-2}}{N^2} + O \left(\frac{1}{N^3}\right).
\end{align}
\end{lemma}

\begin{proof}[Proof of Lemma \ref{lem:gcw_partial}]
Expanding the definitions, we have
\begin{align}
    \partial_i (m(x))^k &= m(x^{+i})^k - m(x^{-i})^k \\
    &= \left( m(x^{-i}) + \frac{1}{N} \right)^k - m(x^{-i})^k \\
    &= \frac{k m(x^{-i})^{k-1}}{N} + O \left(\frac{1}{N^2}\right),
\end{align}
which proves the first assertion.
Similarly, we have
\begin{align}
    \partial_j \partial_i (m(x))^k &= m(x^{+i+j})^k - m(x^{-i+j})^k - m(x^{+i-j})^k + m(x^{-i-j})^k \\
    &= \left( m(x^{-i-j}) + \frac{2}{N} \right)^k - 2 \left( m(x^{-i-j}) + \frac{1}{N} \right)^k + m(x^{-i-j})^k \\
    &= \left( 4 \binom{k}{2} - 2 \binom{k}{2} \right) \frac{m(x^{-i-j})^{k-2}}{N^2} + O \left(\frac{1}{N^3}\right),
\end{align}
proving the second assertion.
\end{proof}

With this lemma in hand we may prove the following contraction estimate for conditioned GCWM Glauber dynamics chains started within
the slightly smaller set $\Lambda_{m_*}^\eps$ for any $m_* \in U_\beta^\GCW$.
Here is where we make use of the strict concavity in the definition of $U_\beta^\GCW$.

\begin{proposition}
\label{prop:gcw_contraction}
For small enough $\eps > 0$,
if $x, y \in \Lambda_{m_*}^\eps$ differ in exactly one coordinate, then
under the monotone coupling we have
\begin{equation}
    \E \left[
        \dh(X_1^x, X_1^y)
    \right] \leq 1 - \frac{\kappa}{N}
\end{equation}
for some constant $\kappa>0$ only depending on $\eps$.
\end{proposition}

\begin{proof}[Proof of Proposition \ref{prop:gcw_contraction}]
Let us denote by $j$ the coordinate at which $x$ and $y$ differ.
Suppose without loss of generality that $x(j) = 0$ and $y(j) = 1$ so that $y = x^{+j}$ and $x = x^{-j}$.
If the coordinate $j$ is chosen to be updated, this discrepancy will be resolved and the resulting Hamming distance will be $0$.
On the other hand, if $i \neq j$ is chosen then the Hamming distance will either remain $1$, or it will increase to $2$ with probability
\begin{equation}
    \left| \varphi(N \partial_i \hamilton(y)) - \varphi(N \partial_i \hamilton(x)) \right|,
\end{equation}
since the new spins at coordinate $i$ are optimally coupled Bernoulli random variables with these means.
Now, applying a second-order Taylor expansion and using the fact that $y = x^{+j}$ and $x = x^{-j}$ we obtain
\begin{align}
    \E \left[ \dh(X_1^x, X_1^y) \right] &= \frac{1}{N} \sum_{i \neq j} \left(
        1 + \left|
            \varphi(N \partial_i \hamilton(y))
            -\varphi(N \partial_i \hamilton(x))
        \right|
    \right) \\
    &= 1 - \frac{1}{N} + \frac{1}{N} \sum_{i \neq j} \left|
        \varphi'(N \partial_i \hamilton(x)) \cdot N \partial_j \partial_i \hamilton(x)
        + O \left( \left( N \partial_j \partial_i \hamilton(x) \right)^2 \right)
    \right|.
\label{eq:dhamint}
\end{align}
Now by Lemma \ref{lem:gcw_partial} since $\hamilton(x) = h(m(x))$, and $x = x^{-j}$, we have
\begin{align}
    N \partial_i \hamilton(x) &= h'(m(x^{-i})) + O(N^{-1}), \\
    N \partial_j \partial_i \hamilton(x) &= \frac{h''(m(x^{-i}))}{N} + O \left( \frac{1}{N^2} \right).
\end{align}
Additionally, since $\varphi'$, $h'$ and $h''$ are differentiable and $|m(x^{-i}) - m_*| < \eps + N^{-1}$ because $x \in \Lambda_{m_*}^\eps$,
we have
\begin{align}
    \varphi'(N \partial_i \hamilton(x)) &= \varphi'(h'(m_*)) + O(\eps + N^{-1}), \\
    N \partial_j \partial_i \hamilton(x) &= \frac{h''(m_*)}{N} + O \left(\frac{\eps}{N} + \frac{1}{N^2} \right),
\end{align}
where the error terms here do not depend on the coordinate $i$.
Let us now plug these approximations into \eqref{eq:dhamint} and use the fact that $\varphi$ is increasing and
$h$ is convex to remove the absolute value sign.
Also using the fact that $\frac{N-1}{N} \leq 1$, we find that
\begin{align}
    \E[\dh(X_1^x,X_1^y)]
    &\leq 1 - \frac{1}{N} +
        \frac{\varphi'(h'(m_*)) h''(m_*)}{N}
    + O \left( \frac{\eps}{N} + \frac{1}{N^2} \right) \\
    &= 1 - \frac{1}{N}  \left(1 - \left. \frac{d}{dm} \varphi(h'(m)) \right|_{m=m_*} \right) + O \left( \frac{\eps}{N} + \frac{1}{N^2} \right).
\end{align}
Now by Lemma \ref{lem:equiv}, since $m_* \in U_\beta^\GCW$, the derivative in the expression above is strictly less than $1$.
Thus the desired conclusion follows if we take $\eps$ small enough.
\end{proof}
\subsection{Proof of Theorem \ref{thm:main_gcwm}}
\label{sec:proof}

In this section we finish the proof of Theorem \ref{thm:main_gcwm}.
Again, for notational convenience, we set $\mu = \calM_{\beta,m_*}^\GCW$ for some fixed $m_* \in M_\beta^\GCW$.
We will apply Theorem \ref{thm:general} with the set $\Lambda = \Lambda_{m_*}^{\eps}$ for small enough $\eps > 0$.
Note that $\mu(\Lambda^\comp) \leq e^{-\Omega(N)}$ by \eqref{eq:lambdalarge_gcw}.
In addition, since we may move between any pair of elements in $\Lambda$ by flipping one bits at a time without
leaving $\Lambda$ as long as $N$ is large enough, we have $\diam{\Lambda} \leq N$.
By Proposition \ref{prop:gcw_contraction}, we may also take $\alpha = 1 - \frac{\kappa}{N}$ in condition \ref{cond:contraction}.
Let us take $T = N^2$ which is less than $\frac{1}{\mu(\Lambda^\comp)} = e^{\Omega(N)}$ for large enough $N$.
Assuming that condition \ref{cond:localFKG} holds for now, the conclusion of Theorem \ref{thm:general} is that for any increasing
$f,g$, we have
\begin{equation}
    \Cov_\mu[f,g] \geq - \| f \|_\infty \| g \|_\infty \cdot 300 N \left(
        N^2 e^{-\Omega(N)}+ \left(1 - \frac{\kappa}{N} + \frac{20}{N^2}\right)^N \cdot N
    \right).
\end{equation}
Now since $1 - \frac{\kappa}{N} + \frac{20}{N^2} \leq \Exp{- \frac{\kappa}{2N}}$ for large enough $N$, the entire expression in parentheses
above is $e^{-\Omega(N)}$, which finishes the proof.

It just remains to check condition \ref{cond:localFKG}, the local FKG condition near $\Lambda$.
Since $\Lambda = \Lambda_{m_*}^\eps$ and $\mu$ is the GCWM measure conditioned on $\Lambda_{m_*}^\eta$ for some $\eta > \eps$,
any $z$ with $\dh(z,\Lambda) \leq 2$ is also in $\Lambda_{m_*}^\eta$ for large enough $N$.
So for any $x,y$ with $\dh(x,\Lambda), \dh(y,\Lambda) \leq 1$, if we also have $\dh(x \vee y, \{x,y\}), \dh(x \wedge y, \{x,y\}) \leq 1$,
then all four states $x,y, x \vee y$ and $x \wedge y$ are in $\Lambda_{m_*}^\eta$, meaning that for $z \in \{x,y,x \vee y, x \wedge y\}$
we have $\mu(z) \propto \Exp{N \hamilton(z)}$, with the same proportionality constant, say $\frac{1}{\calZ}$.
So we have
\begin{align}
    \mu(x \vee y) \mu(x \wedge y) &= \frac{1}{\calZ^2} \Exp{N \hamilton(x \vee y) + N \hamilton (x \wedge y)} , \\
    \mu(x) \mu(y) &= \frac{1}{\calZ^2} \Exp{N \hamilton(x) + N \hamilton (y)}.
\end{align}
Recall that $\hamilton(z) = h(m(z)) = \sum_{j=1}^K \beta_j (m(z))^j$ with $\beta_j \geq 0$ for $j \geq 2$.
For the $j=1$ term we have
\begin{equation}
    m(x \vee y) + m(x \wedge y) = m(x) + m(y)
\end{equation}
by the inclusion-exclusion principle.
Since $\beta_j \geq 0$ for $j \geq 2$, to prove the local FKG condition \ref{cond:localFKG} it thus suffices to show that
\begin{equation}
\label{eq:fkggoal}
    (m(x \vee y))^k + (m(x \wedge y))^k \geq m(x)^k + m(y)^k
\end{equation}
for any $k \geq 2$.
Let us define
\begin{align}
    A &= \frac{1}{N} \left| \{ i : x(i) = 1, y(i) = 1 \} \right|, \\
    B &= \frac{1}{N} \left| \{ i : x(i) = 1, y(i) = 0 \} \right|, \\
    C &= \frac{1}{N} \left| \{ i : x(i) = 0, y(i) = 1 \} \right|.
\end{align}
Then we have $m(x) = A + B$, $m(y) = A + C$, $m(x \vee y) = A + B + C$, and $m(x \wedge y) = A$.
So \eqref{eq:fkggoal} is the same as
\begin{equation}
    (A + B + C)^k + A^k \geq (A+B)^k + (A+C)^k,
\end{equation}
which holds for all $k \geq 1$ and all nonnegative real numbers $A,B,C$ 
by convexity of the function $s \mapsto s^k$.
This proves \eqref{eq:fkggoal} and thus finishes the proof of Theorem \ref{thm:main_gcwm}.
\section{Application to exponential random graph models}
\label{sec:ergm}

In this section we will prove Theorem \ref{thm:main_ergm} by an application of Theorem \ref{thm:general} to exponential random graph models (ERGMs).
In contrast to the generalized Curie--Weiss models studied in Section \ref{sec:cw}, for ERGMs the mechanisms for mixing are somewhat decoupled
from the large deviations principle defining the phase measures.
Nevertheless, recent work has elucidated the relationship between these two concepts enough that we may apply Theorem \ref{thm:general} in a relatively
straightforward manner.
We will introduce the large deviations principle of \cite{chatterjee2013estimating} properly in the next section, after which we will describe
the metastable mixing result from \cite{winstein2025concentration} which we will use as input for Theorem \ref{thm:general}.
This metastable mixing result builds on work done by \cite{bresler2024metastable} which in turn draws from \cite{bhamidi2008mixing}, the first study of the
dynamics of ERGMs.

\subsection{Large deviations principle}
\label{sec:ergm_background_ldp}

Let us first properly introduce the large deviations principle alluded to in Section
\ref{sec:intro_setup_ergm}.
For this, we will need to quickly define graphons and the cut distance; a more complete
picture of this subject can be found in \cite{lovasz2012large}.
Briefly, a graphon is a symmetric measurable function $W : [0,1]^2 \to [0,1]$, which is a limiting
version of an adjacency matrix.
Indeed, representing a graph $x$ as a collection of edge indicators $x(e)$ for $e \in \{ 0, 1 \}^{\edgeset}$, we may
define a graphon corresponding to $x$ as
\begin{equation}
    W_x = \sum_{\{u,v\} \in \edgeset} x(\{u,v\}) \left(
    \indd{\left[ \frac{u-1}{n}, \frac{u}{n} \right) \times \left[ \frac{v-1}{n}, \frac{v}{n} \right)}
    + \indd{\left[ \frac{v-1}{n}, \frac{v}{n} \right) \times \left[ \frac{u-1}{n}, \frac{u}{n} \right)}
    \right).
\end{equation}
We consider graphons to be identical if they only differ on a set of Lebesgue measure zero; thus it is impossible
in general to recover $x$ from $W_x$ unless one also knows how many vertices $x$ has.
Nevertheless, we will often treat a graph $x$ as the graphon $W_x$.
Much like we consider graphs to be isomorphic when we relabel their vertices, we also consider two graphons $W$ and $W'$
to be isomorphic if there is a measure-preserving bijection $\sigma : [0,1] \to [0,1]$ such that
$W'(s,t) = W(\sigma(s),\sigma(t))$ for Lebesgue-almost-every $(s,t) \in [0,1]^2$.

The natural distance on the space of graphons is the \emph{cut distance}
\begin{equation}
    \db(W,W') = \inf_{\substack{\sigma : [0,1] \to [0,1] \\ \text{measure-preserving} \\ \text{bijection}}}
    \sup_{S, T \subseteq [0,1]} \left| \int_S \int_T (W'(s,t) - W(\sigma(s), \sigma(t))) \,\,dt \,ds \right|.
\end{equation}
Note that the absolute value is \emph{outside} of the integral, which allows for cancellations on large enough sets $S$ and $T$,
meaning that, for instance, if $X$ is sampled from the Erd\H{o}s-R\'enyi model $\calG(n,p)$ then $\db(X,W_p)$ is typically small,
where $W_p$ is the constant graphon $W_p(s,t) = p$ and we abuse notation by writing $X$ instead of $W_X$ as defined above.

We may define subgraph densities $t(G,W)$ for a graphon $W$ and a finite graph $G$, with vertex and edge set $\calV(G)$ and $\calE(G)$ respectively,
as follows:
\begin{equation}
    t(G,W) = \int_{(s_v) \in [0,1]^{\calV(G)}} \prod_{\{u,v\} \in \calE(G)} W(s_u, s_v) \prod_{v \in \calV(G)} ds_v.
\end{equation}
Note that this extends the notion of subgraph densities defined in Section \ref{sec:intro_setup_ergm} above, as $t(G,x) = t(G,W_x)$
for any finite graph $x$.
In addition, convergence in cut distance is equivalent to convergence of all subgraph densities:
\begin{equation}
    \db(W_n, W) \to 0 \qquad \Longleftrightarrow \qquad t(G,W_n) \to t(G,W) \text{ for all } G.
\end{equation}
Having introduced the cut distance and given a bit of context, we may now state the large deviations principle of \cite{chatterjee2013estimating}.
Recall that $M_\beta^\ERG \subseteq [0,1]$ denotes the set of global maximizers of the function $L_\beta^\ERG$ defined in \eqref{eq:Lerg_def}.

\begin{theorem}[Theorems 3.2 and 4.1 of \cite{chatterjee2013estimating}]
\label{thm:ergm_ldp}
For any $\eta > 0$, there are constants $c(\eta), C(\eta) > 0$ such that
\begin{equation}
    \P \left[ \inf_{p_* \in M_\beta^\ERG} \db(X,W_{p_*}) > \eta \right] \leq C(\eta) \Exp{- c(\eta) \cdot n^2}.
\end{equation}
\end{theorem}

With this, we may properly define the phase measures discussed in Section \ref{sec:intro_setup_ergm}.
We fix $\eta > 0$ small enough so that the cut distance balls
\begin{equation}
    \ball \coloneqq \left\{ W : \db(W, W_{p_*}) \leq \eta \right\}
\end{equation}
are disjoint as $p_*$ varies in $M_\beta^\ERG$, and for each $p_* \in M_\beta^\ERG$ we define the $n$-vertex phase measure $\calM^\ERG_{\beta,p_*}$
as the full $n$-vertex ERGM measure \eqref{eq:full_ergm_measure} conditioned on $\ball$, or more precisely, on the set
\begin{equation}
    \left\{ x \in \{0,1\}^{\edgeset} : W_x \in \ball \right\}.
\end{equation}
This is the measure to which Theorem \ref{thm:main_ergm} applies.
\subsection{Contraction in a large connected set}
\label{sec:ergm_background_mixing}

Let us now introduce the contraction property which will be crucial for applying Theorem \ref{thm:general}
to exponential random graph models.
This originally derives from work of \cite{bhamidi2008mixing}, which proves rapid mixing of the ERGM Glauber dynamics
in the \emph{subcritical} or \emph{high temperature} regime of parameters (which we avoid defining as it will not be relevant for us).
In brief, \cite{bhamidi2008mixing} shows that the dynamics exhibits contraction in a ``good'' set $\Gamma_{p_*}^\eps$ (to be defined below)
which is also a well of attraction, leading to rapid mixing.
Later, \cite{bresler2024metastable} showed that $\Gamma_{p_*}^\eps$ is also \emph{large} under the phase measure $\calM^\ERG_{\beta,p_*}$
for any $p_* \in U_\beta^\ERG$,
the set of global maximizers of the function $L_\beta^\ERG$ (defined in \eqref{eq:Lerg_def}) for which $L_\beta^\ERG$ is strictly concave.
This led to \emph{metastable mixing} where the dynamics mixes rapidly to stationarity from a \emph{warm start} with $X_0 \sim \calG(n,p_*)$,
the Erd\H{o}s--R\'enyi measure.
This situation was further analyzed by \cite{winstein2025concentration} to extract a \emph{connected} subset $\Lambda$ of $\Gamma_{p_*}^\eps$
which is still large.
Connectedness was an important property in the use of path coupling starting from a not necessarily warm start, which useful in
\cite{winstein2025concentration} for deriving a concentration inequality for the phase measures, along the lines of \cite{ganguly2024sub}
which derives such an inequality for the subcritical regime of parameters.
As we have seen, connectedness is also important for the use of our Theorem \ref{thm:general}, for much the same reason.

To be precise and avoid confusion, in this section we use $(Y_t^x)$ to denote the unconditioned ERGM Glauber dynamics started at $x$
(i.e.\ the Glauber dynamics with respect to the full ERGM measure \eqref{eq:full_ergm_measure}).
This only differs from the conditioned dynamics $(X_t^x)$ when it would attempt to leave the cut distance ball $\ball$, in which case the unconditioned
dynamics proceeds and the conditioned dynamics stays put.

To define the connected set $\Lambda$ let us first recall the ``good'' set $\Gamma_{p_*}^\eps$ defined in \cite{bhamidi2008mixing}.
Intuitively, this is the set where the ERGM Glauber dynamics behaves much like simply resampling edges with probability $p_*$ of inclusion.
To make this precise, for any finite graph $G = (\calV,\calE)$ and any edge $e \in \edgeset$, let us define
\begin{equation}
\label{eq:rgxe}
    r_G(x, e) \coloneqq \left( \frac{\partial_e N_G(x)}{2 |\calE| n^{|\calV|-2}} \right)^{\frac{1}{|\calE|-1}},
\end{equation}
recalling the definition of the homomorphism count $N_G(x)$ from Section \ref{sec:intro_setup_ergm}
and the discrete derivative from Section \ref{sec:cw_contraction}.
To motivate the definition \eqref{eq:rgxe}, note that if $X \sim \calG(n,p)$ then $r_G(X,e) \approx p$ with high probability.
This is because $\partial_e N_G(x)$ is the number of homomorphisms of $G$ in $x^{+e}$ which make use of the edge $e$;
there are $2 |\calE|$ edges in $G$ which can map to $e$ (counting both orientations), and the remaining $|\calV|-2$ vertices
must map somewhere in $[n]$.
In $X \sim \calG(n,p)$, the probability that any such map is a homomorphism in $X^{+e}$ is thus $p^{|\calE|-1}$ since
one edge is already included automatically.
This shows that $r_G(X,e) \approx p$ as stated.

Additionally, by a derivation much like the one laid out for GCWMs at the beginning of Section \ref{sec:cw_contraction},
the update probabilities of the unconditioned ERGM Glauber dynamics take the form
\begin{equation}
    \P \left[ Y_1^x(e) = 1 \right] = \varphi(n^2 \partial_e \hamilton_\beta^\ERG(x))
    = \varphi \left( \sum_{j=0}^K 2 \beta_j |\calE_j| \cdot r_{G_j}(x,e)^{|\calE_j|-1} \right),
\end{equation}
recalling the definition of the Hamiltonian \eqref{eq:ERGM_hamiltonian},
where we recall the definition of $\varphi(s) = \frac{e^s}{1+e^s}$ from Section \ref{sec:cw} above,
and $G_j = (\calV_j, \calE_j)$ are the graphs in the ERGM specification for $0 \leq j \leq K$.
As in the case of Lemma \ref{lem:equiv} for GCWMs,
one may easily check that each $p_* \in U_\beta^\ERG$ is an attracting fixed point of the map
\begin{equation}
    p \mapsto \varphi(2 h'(p)) = \varphi \left( \sum_{j=0}^K 2 \beta_j |\calE_j| p^{|\calE_j|-1}\right),
\end{equation}
where we write $h = h_\beta^\ERG$ as defined in Section \ref{sec:intro_setup_ergm} (the factor of $2$ arises from the
factor of $\frac{1}{2}$ on the entropy term of $L_\beta^\ERG$).
So if all $r_G(x,e)$ are close to some $p_* \in U_\beta^\ERG$ then, heuristically speaking, $Y_1^x$ will also have this property.
Thus, for any $p_* \in U_\beta^\ERG$ and any $\eps > 0$ we define
\begin{equation}
    \Gamma_{p_*}^\eps \coloneqq \left\{
        x \in \{0,1\}^{\edgeset} :
        \begin{array}{c}
        \left| r_G(x,e) - p_* \right| \leq \eps \text{ for all } e \in \edgeset \text{ and for all}\\
        \text{graphs } G \text{ with } |\calV(G)| \leq \max \{ |\calV_j| : 0 \leq j \leq K \}
        \end{array}
    \right\}.
\end{equation}
With this definition, using a strategy similar to the one presented in Proposition \ref{prop:gcw_contraction} above, one can
show a one-step contraction estimate for the (unconditioned) ERGM Glauber dynamics starting in $\Gamma_{p_*}^\eps$.
The following holds for any $p_* \in U_\beta^\ERG$, and in fact for any strictly concave \emph{local} maximizer of $L_\beta^\ERG$,
although that will not be relevant for us.

\begin{lemma}[Lemma 18 of \cite{bhamidi2008mixing}]
\label{lem:ergm_contraction}
For all small enough $\eps > 0$,
if $x \in \Gamma_{p_*}^\eps$ then for all $e \in \edgeset$, under the monotone coupling, we have
\begin{equation}
    \E \left[ \dh \left( Y_1^{x^{+e}}, Y_1^{x^{-e}} \right) \right] \leq 1 - \frac{\kappa}{n^2},
\end{equation}
for some constant $\kappa > 0$ depending only on the ERGM specification.
\end{lemma}

As mentioned above, it was shown in \cite{bresler2024metastable} that $\Gamma_{p_*}^\eps$ is also large under the phase measure $\calM^\ERG_{\beta,p_*}$
for any $p_* \in U_\beta^\ERG$, leading to metastable mixing.
The following result refines this to show that there is in fact a large \emph{connected} subset of $\Gamma_{p_*}^\eps$.
For notational convenience, let us denote $\calM^\ERG_{\beta,p_*}$ simply by $\mu$ in what follows.

\begin{proposition}[Proposition 5.9 of \cite{winstein2025concentration}]
\label{prop:ergm_goodset}
For all small enough $\eps > 0$ there is a connected set $\Lambda \subseteq \halfball \cap \Gamma_{p_*}^\eps$ with $\diam{\Lambda} \leq 2 n^2$
and $\mu(\Lambda^\comp) \leq e^{-\Omega(n)}$.
\end{proposition}

The importance of the fact that $\Lambda \subseteq \halfball$ in addition to $\Gamma_\eps$ is that when the dynamics starts in $\halfball$,
it will have no chance to attempt to leave $\ball$, at least for large enough $n$.
This means that if $\dh(x, \Lambda) \leq 1$ then $Y_1^x$ and $X_1^x$ (one step of the unconditioned and conditioned dynamics respectively)
have exactly the same distribution, so we may use the contraction estimate given by Lemma \ref{lem:ergm_contraction} for the conditioned
dynamics as well.
With the above preparation, we are ready to prove Theorem \ref{thm:main_ergm}.
\subsection{Proof of Theorem \ref{thm:main_ergm}}
\label{sec:ergm_proof}

Again, for notational convenience, let us set $\mu = \calM^\ERG_{\beta,p_*}$ for some fixed $p_* \in U_\beta^\ERG$,
and let us set $\hamilton = \hamilton_\beta^\ERG$.
We will apply Theorem \ref{thm:general} with the set $\Lambda$ as given by Proposition \ref{prop:ergm_goodset}, with $\diam{\Lambda} \leq 2 n^2$
and $\mu(\Lambda^\comp) \leq e^{-\Omega(n)}$.
Since $\Lambda \subseteq \Gamma_{p_*}^\eps$, Lemma \ref{lem:ergm_contraction} implies that we may take $\alpha = 1 - \frac{\kappa}{n^2}$
in condition \ref{cond:contraction}.
Let us also set $T = n^3$, which is less than $\frac{1}{\mu(\Lambda^\comp)} = e^{\Omega(n)}$ for large enough $n$.
Assuming that condition \ref{cond:localFKG} holds for now, the conclusion of Theorem \ref{thm:general} is then that for any increasing $f,g$,
we have
\begin{equation}
    \Cov_\mu[f,g] \geq - \| f \|_\infty \| g \|_\infty \cdot 300 n^{1.5}\left(
        n^3 e^{-\Omega(n)} + \left( 1 - \frac{\kappa}{n^2} + \frac{20}{n^3} \right)^{n^3} 2 n^2
    \right).
\end{equation}
Now, since $1 - \frac{\kappa}{n^2} + \frac{20}{n^3} \leq \Exp{- \frac{\kappa}{2 n^2}}$ for large enough $n$, the entire expression in parentheses is $e^{-\Omega(n)}$,
finishing the proof.

It just remains to check condition \ref{cond:localFKG}, the local FKG condition near $\Lambda$.
Since $\Lambda \subseteq \halfball$, any $z \in \{0,1\}^{\edgeset}$ with $\dh(z,\Lambda) \leq 2$ is also in $\ball$ for large enough $n$.
Thus for any $x,y$ with $\dh(x,\Lambda), \dh(y,\Lambda) \leq 1$, if we also have $\dh(x \vee y, \{x,y\}), \dh(x \wedge y, \{x,y\}) \leq 1$, then all four states
$x, y, x \vee y$, and $x \wedge y$ are in $\ball$.
Thus for $z \in \{x,y, x \vee y, x \wedge y\}$, we have $\mu(z) \propto \Exp{n^2 \hamilton(z)}$, with the same proportionality constant, say $\frac{1}{\mathcal{Z}}$.
So we find that
\begin{align}
    \mu(x \vee y) \mu (x \wedge y) &= \frac{1}{\mathcal{Z}^2} \Exp{n^2 \hamilton(x \vee y) + n^2 \hamilton(x \wedge y)},  \\
    \mu(x) \mu (y) &= \frac{1}{\mathcal{Z}^2} \Exp{n^2 \hamilton(x) + n^2 \hamilton(y)}.
\end{align}
Recalling the definition from \eqref{eq:ERGM_hamiltonian}, the Hamiltonian can be written as
\begin{equation}
    \hamilton(z) = \sum_{j=0}^K \beta_j \cdot t(G_j,z)
\end{equation}
with $\beta_j \geq 0$ for $j \geq 1$ and $\beta_0 \in \R$, where $t(G,z) = \frac{N_G(z)}{n^{|\calV(G)|}}$ and $N_G(z)$ is the count of homomorphisms
of $G$ in $z$.
Since $G_0$ is a single edge, $N_{G_0}(z)$ is just two times the total edge count in $z$.
So by inclusion-exclusion we have
\begin{equation}
    N_{G_0}(x \vee y) + N_{G_0}(x, \wedge y) = N_{G_0}(x) + N_{G_0}(y),
\end{equation}
and so to verify condition \ref{cond:localFKG}, by positivity of $\beta_j$ for $j \geq 1$, it only remains to show that
\begin{equation}
\label{eq:incexagain}
    N_{G}(x \vee y) + N_{G}(x, \wedge y) \geq N_{G}(x) + N_{G}(y),
\end{equation}
for any graph $G$.
For this, note that for each map $\calV(G) \to [n]$, if the map is a homomorphism in $x$ and $y$ then it is a homomorphism in $x \wedge y$,
and if it is a homomorphism in $x$ or in $y$ then it is a homomorphism in $x \vee y$.
Thus \eqref{eq:incexagain} follows from the inclusion-exclusion principle again.
Note that there may be homomorphisms in $x \vee y$ which are not homomorphisms in either $x$ or $y$, which is why the above is an inequality, not an equality.
In any case, this validates condition \ref{cond:localFKG} and thus finishes the proof of Theorem \ref{thm:main_ergm}.
\section{Central limit theorems for phase-coexistence ERGMs}
\label{sec:covariance}

In this section we state a useful consequence of nearly-positive correlations of increasing functions and use it, along with the work of
\cite{winstein2025concentration,winstein2025quantitative}, to prove quantitative central limit theorems for ERGMs in the \emph{phase coexistence}
regime in general, without the so-called forest assumption that was made in \cite{winstein2025quantitative}.
This answers a question posed by \cite{bianchi2024limit} and completes the picture of central limit theorems for ferromagnetic ERGMs,
aside from the \emph{critical} case where we do not expect such central limit
theorems to hold (for instance, see \cite{mukherjee2013statistics}).
The main result that allows us to complete this program is an extension of a classical covariance inequality discussed in the next section,
which we expect to be useful in a variety of other contexts as well.

\subsection{Extension of a classical covariance inequality}
\label{sec:covariance_extension}

In this section we prove a theorem which allows us to bound covariances of general functions by covariances of individual spins, under measures satisfying an approximate FKG inequality.
The exact version of this (assuming the non-approximate FKG inequality) was first mentioned by \cite{newman1980normal} without proof, and with an extraneous constant factor of $3 \sqrt{2}$.
A proof was given by \cite{bulinski1998asymptotical} in Russian without this constant factor, later translated to English in the book \cite{bulinski2007limit}.
Our proof of the approximate version mirrors the proof from \cite{bulinski2007limit}.
In what follows, we consider measures $\mu$ supported on the $N$-fold product of a particular \emph{finite interval} $I \subseteq \R$, and we denote by $|I|$ the length of the interval.
In addition, for a function $F : I^N \to \R$ we define
\begin{equation}
    \Lip_i(F) = \sup \left\{ \frac{|F(x) - F(y)|}{|x(i)-y(i)|} : x,y \in I^N, x(i) \neq y(i), x(j) = y(j) \text{ for } j \neq i \right\}.
\end{equation}
With this definition, we may state the covariance bound.

\begin{theorem}
\label{thm:cov}
Suppose that there is some $\delta \geq 0$ such that the measure $\mu$ on $I^N$ satisfies
\begin{equation}
    \Cov_\mu[f,g] \geq - \| f \|_\infty \| g \|_\infty \delta
\end{equation}
for any coordinate-wise increasing $f,g$.
Then for any functions $F,G : I^N \to \R$ (not necessarily increasing), we have
\begin{equation}
    \left| \Cov_\mu[F,G] \right| \leq \sum_{i,j=1}^N \Lip_i(F) \Lip_j(G) \left( \Cov_\mu[X(i),X(j)] + 4|I|^2 \delta \right),
\end{equation}
where $X \in I^N$ denotes a sample from $\mu$.
\end{theorem}

\begin{proof}[Proof of Theorem \ref{thm:cov}]
Let us define auxiliary functions $F_\pm, G_\pm : I^N \to R$ by
\begin{equation}
    F_\pm (x) = F(x) \pm \sum_{i=1}^N \Lip_i(F) \cdot x(i)
    \qquad \text{and} \qquad
    G_\pm (x) = G(x) \pm \sum_{j=1}^N \Lip_j(G) \cdot x(j).
\end{equation}
Then $F_+$  and $G_+$ are coordinate-wise increasing and $F_-$ and $G_-$ are coordinate-wise decreasing.
So by the hypothesis of the approximate FKG inequality we have
\begin{align}
    - \| F_+ \|_\infty \| G_+\|_\infty \delta &\leq \Cov_\mu[F_+, G_+] \notag \\
    &= \Cov_\mu[F, G] & \qquad \eqcolon A \\
    &\qquad + \sum_{i,j=1}^N \Lip_i(F) \Lip_j(G) \Cov_\mu[X(i), X(j)] & \qquad \eqcolon B \\
    &\qquad + \sum_{i=1}^N \Lip_i(F) \Cov_\mu[X(i), G(X)] & \qquad \eqcolon C \\
    &\qquad + \sum_{j=1}^N \Lip_j(G) \Cov_\mu[F(X), X(j)]. & \qquad \eqcolon D
\end{align}
Let us use the letters $A,B,C,D$ to denote the expressions on the lines above corresponding to their labels.
So the above inequality is
\begin{equation}
\label{eq:cov1}
    A + B + C + D \geq - \| F_+ \|_\infty \| G_+ \|_\infty \delta.
\end{equation}
Then similar algebra (using the fact that $-F_-$ and $- G_-$ are coordinate-wise increasing) shows that
\begin{align}
\label{eq:cov2}
    A + B - C - D &= \Cov_\mu[- F_-, - G_-] \geq - \| F_- \|_\infty \| G_- \|_\infty \delta, \\
\label{eq:cov3}
    - A + B - C + D &= \Cov_\mu[F_+, -G_-] \geq - \| F_+ \|_\infty \| G_- \|_\infty \delta, \\
\label{eq:cov4}
    - A + B + C - D &= \Cov_\mu[-F_-, G_+] \geq - \| F_- \|_\infty \| G_+ \|_\infty \delta.
\end{align}
Adding \eqref{eq:cov1} and \eqref{eq:cov2} we arrive at
\begin{equation}
\label{eq:finalcov1}
    -A \leq B + \frac{1}{2} \left( \| F_+ \|_\infty \| G_+ \|_\infty + \| F_- \|_\infty \| G_- \|_\infty \right) \delta,
\end{equation}
and adding \eqref{eq:cov3} and \eqref{eq:cov4} we arrive at
\begin{equation}
\label{eq:finalcov2}
    A \leq B + \frac{1}{2} \left( \| F_+ \|_\infty \| G_- \|_\infty + \| F_- \|_\infty \| G_+ \|_\infty \right) \delta.
\end{equation}
Now note that we may shift $F_\pm$, and $G_\pm$ by any constants (which may be different for each function) without changing the above derivation,
since all of the above is only in terms of covariances.
So, if $I = [a,b]$ and we shift each function so that it takes the value $0$ at $(a,\dotsc,a)$, then we obtain \eqref{eq:finalcov1} and \eqref{eq:finalcov2}
with $F_\pm$ and $G_\pm$ replaced by $\widetilde{F_\pm}$ and $\widetilde{G_\pm}$ which satisfy
\begin{equation}
    \| \widetilde{F_\pm} \|_\infty \leq \sum_{i=1}^N |I| 2\Lip_i(F),
    \qquad \text{and} \qquad
    \| \widetilde{G_\pm} \|_\infty \leq \sum_{j=1}^N |I| 2\Lip_i(G).
\end{equation}
Indeed, we have $\Lip_i(F_\pm) \leq 2 \Lip_i(F)$ and similarly for $G$, and any element $x \in [a,b]^N$
can be reached by traveling along coordinates one-by-one and incurring the corresponding increase which is bounded by the Lipschitz constant.
Plugging these bounds into \eqref{eq:finalcov1} and \eqref{eq:finalcov2} finishes the proof.
\end{proof}
\subsection{Finishing the proof of ERGM central limit theorems}
\label{sec:covariance_clt}

Let us now turn to the proof of Theorem \ref{thm:main_clt}, which answers a question of \cite{bianchi2024limit}
by giving a central limit theorem for the edge count of an $n$-vertex sample $X \sim \calM^\ERG_{\beta,p_*}$ for any $p_* \in U_\beta^\ERG$
(as defined in Section \ref{sec:intro_setup_ergm}).
For notational convenience, we let $X$ denote such a sample throughout the rest of this section.
Recall that Theorem \ref{thm:main_clt} states that
\begin{equation}
    \distance \left( \frac{|\calE(X)| - \E[|\calE(X)|]}{\sigma_n}, Z \right) \lesssim n^{- \frac{1}{2} + \eps},
\end{equation}
where $|\calE(X)|$ is the edge count, $\sigma_n$ is an explicit quantity to be defined below,
and $\distance$ denotes either the 1-Wasserstein distance $\dW$ or the Kolmogorov distance $\dK$
between real-valued random variables, which are defined as follows:
\begin{equation}
    \dW(S,Z) \coloneqq \inf_{\substack{S' \stackrel{d}{=} S \\ Z' \stackrel{d}{=} Z}} \E \left[ | S' - Z' |\right]
    \qquad \text{and} \qquad
    \dK(S,Z) \coloneqq \sup_{s \in \bbR} \left| \P[S \leq s] - \P[Z \leq s] \right|.
\end{equation}
The infimum in the definition of $\dW$ is taken over all couplings $(S',Z')$ of variables with the same marginal distributions
as $S$ and $Z$.

Let us write down the formula for the quantity $\sigma_n$ used in the statement of Theorem \ref{thm:main_clt}.
First, write $h = h_\beta^\ERG$ for the polynomial which was defined in Section \ref{sec:intro_setup_ergm}, i.e.
\begin{equation}
    h(p) = h_\beta^\ERG(p) = \sum_{j=0}^K \beta_j p^{|\calE_j|},
\end{equation}
where $G_j = (\calV_j,\calE_j)$ denote the graphs in the specification of the ERGM, recalling that $G_0$ is always a single edge.
Now for a fixed choice of $p_* \in U_\beta^\ERG$, we define
\begin{equation}
\label{eq:sigmadef}
    \sigma_n^2 \coloneqq \frac{p_*(1-p_*) \binom{n}{2}}{1 - 2 p_*(1-p_*) h''(p_*)}.
\end{equation}
We will prove Theorem \ref{thm:main_clt} with this choice of $\sigma_n = \sqrt{\sigma_n^2}$.
Note that $\sigma_n^2$ is strictly larger than the corresponding variance $p_* (1-p_*) \binom{n}{2}$ of the edge count in the Erd\H{o}s--R\'enyi
model $\calG(n,p_*)$.
Our methods will also yield the following corollary for the subgraph homomorphism count $N_G(X)$ defined in Section \ref{sec:intro_setup_ergm}.

\begin{corollary}
\label{cor:clt}
Let $X \sim \calM^\ERG_{\beta,p_*}$ for some $p_* \in U_\beta^\ERG$, and
let $G = (\calV,\calE)$ be a fixed finite graph.
Then for any $\eps > 0$, we have
\begin{equation}
    \dW \left( \frac{N_G(X) - \E[N_G(X)]}{2 |\calE| p_*^{|\calE|-1} n^{|\calV|-2} \sigma_n}, Z \right) \lesssim n^{-\frac{1}{2} + \eps}.
\end{equation}
\end{corollary}

In addition, our methods yield quantitative central limit theorems for the degree of any particular vertex and for the number of subgraph homomorphisms which
use a particular vertex; these are ``local'' versions of Theorem \ref{thm:main_clt} and Corollary \ref{cor:clt}.
To avoid introducing all of the required notation, we do not explicitly state these local results, but they take the exact same form as given in
\cite[Theorem 1.3 and Corollary 1.4]{winstein2025quantitative}, but without the requirement that we are in the \emph{forest or phase-uniqueness}
regime assumed in that work.

As mentioned in \cite[near the proofs of Lemmas 4.2 and 4.5]{winstein2025quantitative}, the missing pieces for the above central limit theorems in the \emph{non-forest}, \emph{phase-coexistence} regime
were only the following two propositions, the first of which gives a bound on the multi-way covariance of edge indicators and the second of which gives a bound
on the distance between the density in the ERGM and the limiting density $p_*$.

\begin{proposition}
\label{prop:multilinear}
Let $X \sim \calM^\ERG_{\beta,p_*}$ for some $p_* \in U_\beta^\ERG$.
For any fixed $k$, let $e_1, \dotsc, e_k \in \edgeset$ be distinct potential edges.
Then
\begin{equation}
    \left|
        \E \left[ \prod_{j=1}^k X(e_j) \right] - \E[X(e)]^k
    \right| \lesssim \frac{1}{n},
\end{equation}
where $e$ is an arbitrary edge in $\edgeset$.
\end{proposition}

\begin{proposition}
\label{prop:marginal}
Let $X \sim \calM^\ERG_{\beta,p_*}$ for some $p_* \in U_\beta^\ERG$.
For any $e \in \edgeset$,
\begin{equation}
    \left| \E[X(e)] - p_* \right| \lesssim \sqrt{\frac{\log n}{n}}.
\end{equation}
\end{proposition}

As was shown in \cite{ganguly2024sub} for the subcritical case, Proposition \ref{prop:marginal} is a consequence of Proposition \ref{prop:multilinear},
in addition to the concentration inequality \cite[Theorem 1]{ganguly2024sub}.
This concentration inequality was already extended to the phase measure $\calM^\ERG_{\beta,p_*}$ in \cite[Theorem 1.1]{winstein2025concentration},
and the derivation of Proposition \ref{prop:marginal} from Proposition \ref{prop:multilinear} for $\calM^\ERG_{\beta,p_*}$ already appears in \cite[Section 6.2]{winstein2025concentration},
so we will simply prove Proposition \ref{prop:multilinear} using Theorem \ref{thm:cov}, and this will finish the proof of Theorem \ref{thm:main_clt} and
Corollary \ref{cor:clt}, as well as both of the local results mentioned above which are extensions of \cite[Theorem 1.3 and Corollary 1.4]{winstein2025quantitative}.

\begin{remark}
Follow up work by the second author \cite{winstein2026wasserstein} improves the
upper bound given in Proposition \ref{prop:marginal} to $\frac{1}{n}$,
using as input our Theorem \ref{thm:main_ergm}, but via a different
line of reasoning.
\end{remark}

In order to prove Proposition \ref{prop:multilinear},
we will make use of the following bound on covariances of individual edge indicators which holds within phase measures as it follows simply from
the aforementioned concentration inequality \cite[Theorem 1.1]{winstein2025concentration}.

\begin{lemma}[Lemma 6.1 of \cite{winstein2025concentration}]
\label{lem:cov}
Let $X \sim \calM^\ERG_{\beta,p_*}$ for some $p_* \in U_\beta^\ERG$.
For any distinct $e, e' \in \edgeset$, we have
\begin{equation}
    \left| \Cov[X(e),X(e')] \right| \lesssim \frac{1}{n}.
\end{equation}
Moreover, if $e$ and $e'$ do not share a vertex, then the right-hand side above can be improved to $\frac{1}{n^2}$.
\end{lemma}

Our proof of Proposition \ref{prop:multilinear} mirrors that of \cite{ganguly2024sub}, but now plugging in
Theorem \ref{thm:cov} instead of the classical result of \cite{newman1980normal,bulinski1998asymptotical},
and using Lemma \ref{lem:cov} to bound the covariances of individual edge variables.
Note also that our proof of Proposition \ref{prop:multilinear} addresses \cite[Conjecture 6.3]{winstein2025concentration}.

\begin{proof}[Proof of Proposition \ref{prop:multilinear}]
Consider the following telescoping sum representation for the quantity in question:
\begin{align}
    \left| \E \left[ \prod_{j=1}^k X(e_j) \right] - \E[X(e)]^k \right|
    &= \left|
        \sum_{j=1}^{k-1} \left( \E \left[ \prod_{i=1}^{j+1} X(e_i) \right] \E[X(e)]^{k-1-j} - \E\left[ \prod_{i=1}^j X(e_i) \right] \E[X(e)]^{k-j} \right)
    \right| \\
    &= \left| \sum_{j=1}^{k-1} \Cov \left[ X(e_1) \dotsb X(e_j), X(e_{j+1}) \right] \cdot \E[X(e)]^{k-1-j} \right| \\
    &\leq \sum_{j=1}^{k-1} \left| \Cov \left[ X(e_1) \dotsb X(e_j), X(e_{j+1}) \right] \right|,
\end{align}
using at the last step that $\E[X(e)] \leq 1$.

Now we will show that each term in the above sum is $\lesssim \frac{1}{n}$, which will finish the proof.
\begin{equation}
    F(x) = x(e_1) \dotsb x(e_j)
    \qquad \text{and} \qquad
    G(x) = x(e_{j+1}).
\end{equation}
Then if $e \in \{e_1,\dotsb,e_j\}$ we have $\Lip_e(F) = 1$ and otherwise we have $\Lip_e(F) = 0$.
Similarly, $\Lip_e(G) = 1$ if and only if $e = e_{j+1}$.
Now by Theorem \ref{thm:main_ergm}, we may take $\delta = e^{-\Omega(n)}$ in the hypothesis of Theorem \ref{thm:cov}, and so we obtain
\begin{equation}
    \left| \Cov[F(X), G(X)] \right| \leq \sum_{i=1}^j \left( \Cov[X(e_i), X(e_{j+1})] + e^{-\Omega(n)} \right) \lesssim \frac{1}{n}
\end{equation}
by Lemma \ref{lem:cov}.
This finishes the proof.
\end{proof}

\bibliographystyle{plain}
\bibliography{biblio}

\end{document}